\documentstyle{article}

\newcommand{\aA}{\mbox{$\cal A$}}

\newcommand{\eL}{\mbox{$\cal L$}}

\newcommand{\eS}{\mbox{$\cal S$}}
\newcommand{\kon}{\wedge}
\newcommand{\dis}{\vee}
\newcommand{\str}{\rightarrow}

\newcommand{\mj}{\mbox{\bf 1}}
\newcommand{\mO}{\mbox{\bf 0}}
\newcommand{\df}{\mbox{\scriptsize{\it df}}}
\newcommand{\HDS}{\vrule width0pt height2.3ex depth1.05ex\displaystyle}

\newcommand{\ml}{\mbox{\bf L}}

\newcommand{\eM}{\mbox{$\cal M$}}

\newcommand{\Bx}{{\raisebox{-1pt}{
\begin{picture}(8,8)
\put(3,2){\makebox(0,0){$\Box$}}
\put(3,3.5){\makebox(0,0){$\times$}}
\end{picture}}}\,}

\newcommand{\Bk}{{\raisebox{-1pt}{
\begin{picture}(8,8)
\put(3,2){\makebox(0,0){$\Box$}}
\put(2.9,3.3){\makebox(0,0){\scriptsize +}}
\end{picture}}}\,}

\newcommand{\bx}{{\raisebox{-1pt}{
\begin{picture}(4,4)
\put(2,2){\makebox(0,0){\scriptsize $\Box$}}
\put(2,3){\makebox(0,0){\tiny $\times$}}
\end{picture}}}}

\newcommand{\bk}{{\raisebox{-1pt}{
\begin{picture}(4,4)
\put(2,2){\makebox(0,0){\scriptsize $\Box$}}
\put(1.75,2.8){\makebox(0,0){\tiny +}}
\end{picture}}}}

\def\cirk{\,{\raisebox{.3ex}{\tiny $\circ$}}\,}
\def\ks{\mbox{\footnotesize$\;\xi\;$}}

\def\f#1#2{{{\HDS #1}\over{\HDS #2}}}

\def\set{\mbox{\it Set}}
\def\prop#1#2{\vspace{2ex} \noindent{\sc #1.} {\it #2} \par \vspace{2ex}}
\def\dkz{\noindent{\sc Proof. }}
\def\qed{\hfill $\dashv$}

\def\pl{\!+\!}
\def\mn{\!-\!}

\def\mlk{\hat{\,\ml}}
\def\mld{\check{\,\ml}}

\def\kst{\raisebox{1pt}{\mbox{\tiny$\xi$}}}

\begin{document}

\title{Bicartesian Coherence Revisited}
\author{\small {\sc Kosta Do\v sen} and {\sc Zoran Petri\' c}
\\[1ex]
{\small Mathematical Institute, SANU}\\[-.5ex]
{\small Knez Mihailova 35, p.f.\ 367, 11001 Belgrade,
Serbia}\\[-.5ex]
{\small email: \{kosta, zpetric\}@mi.sanu.ac.yu}}
\date{}
\maketitle

\begin{abstract}
\noindent A survey is given of results about coherence for
categories with finite products and coproducts. For these results,
which were published previously by the authors in several places,
some formulations and proofs are here corrected, and matters are
updated. The categories investigated in this paper formalize
equality of proofs in classical and intuitionistic
conjunctive-disjunctive logic without distribution of conjunction
over disjunction.

\end{abstract}

\vspace{0.3cm} \noindent {\small \emph{Mathematics Subject
Classification (2000)}: 18A30,
18A15, 03G30, 03G10, 03F05, 03F07, 03B20}
\vspace{0.5cm}

\section{\large\bf Coherence}
Categorists call {\it coherence} what logicians would probably
call {\it completeness}. This is, roughly speaking, the question
whether we have assumed for a particular brand of categories all
the equations between arrows we should have assumed. Completeness
need not be understood here as completeness with respect to
models. We may have also a syntactical notion of
completeness---something like the Post completeness of the
classical propositional calculus---but often some sort of
model-theoretical completeness is implicit in coherence questions.
Matters are made more complicated by the fact that categorists do
not like to talk about syntax, and do not perceive the problem as
being one of finding a match between syntax and semantics. They do
not talk of formal systems, axioms and models.

Moreover, questions that logicians would consider to be questions
of {\it decidability}, which is of course not the same as
completeness, are involved in what categorists call coherence. A
coherence problem often involves the question of deciding whether
two terms designate the same arrow, i.e.\ whether a diagram of
arrows commutes. Coherence is understood mostly as solving this
problem, which we call the \emph{commuting problem}, in
\cite{LS86} (see p.\ 117, which mentions \cite{L68} and \cite{L69}
as the origin of this understanding). The commuting problem seems
to be involved also in the understanding of coherence of
\cite{KL80} (Section 10).

Completeness and decidability, though distinct, are not foreign to
each other. A completeness proof with respect to a manageable
model may provide, more or less immediately, tools to solve
decision problems. For example, the completeness proof for the
classical propositional calculus with respect to the two-element
Boolean algebra provides immediately a decision procedure for
theoremhood.

The simplest coherence questions are those where it is intended
that all arrows of the same type should be equal, i.e.\ where the
category envisaged is a preorder. The oldest coherence problem is
of that kind. This problem has to do with monoidal categories, and
was solved by Mac Lane in \cite{ML63}. The monoidal category
freely generated by a set of objects is a preorder.  So Mac Lane
could claim that showing coherence is showing that ``all diagrams
commute''.

In cases where coherence amounts to showing preorder, i.e.\
showing that from a given set of equations, assumed as axioms, we
can derive all equations (provided the equated terms are of the
same type), from a logical point of view we have to do with {\it
axiomatizability}. We want to show that a decidable set of axioms
(and we wish this set to be as simple as possible, preferably
given by a finite number of axiom schemata) delivers all the
intended equations. If preorder is intended, then all equations
are intended. Axiomatizability is in general connected with
logical questions of completeness, and a standard logical notion
of completeness is completeness of a set of axioms. Where all
diagrams should commute, coherence does not seem to be a question
of model-theoretical completeness, but even in such cases it may
be conceived that the model involved is a discrete category.

Categorists are interested in axiomatizations that permit
extensions. These extensions are in a new language, with new
axioms, and such extensions of the axioms of monoidal categories
need not yield preorders any more. Categorists are also
interested, when they look for axiomatizations, in finding the
combinatorial building blocks of the matter. The axioms are such
building blocks, as in knot theory the Reidemeister moves are the
combinatorial building blocks of knot and link equivalence (see
\cite{BZ85}, Chapter 1, or any other textbook in knot theory).

In Mac Lane's second coherence result of \cite{ML63}, which has to
do with symmetric monoidal categories, it is not intended that all
equations between arrows of the same type should hold. What Mac
Lane does can be described in logical terms in the following
manner. On the one hand, he has an axiomatization, and, on the
other hand, he has a model category where arrows are permutations;
then he shows that his axiomatization is complete with respect to
this model. It is no wonder that his coherence problem reduces to
the completeness problem for the usual axiomatization of symmetric
groups.

Algebraists do not speak of axiomatizations, but of {\it
presentations by generators and relations}. The axiomatizations we
envisage are purely equational axiomatizations, as in algebraic
varieties. Such were the axiomatizations of \cite{ML63}.
Categories are algebras with partial operations, and we are
interested in the equational theories of these algebras.

In Mac Lane's coherence results for monoidal and symmetric
monoidal categories one has to deal only with natural
isomorphisms. However, in the coherence result for symmetric
monoidal closed categories of \cite{KML71} there are already
natural and dinatural transformations that are not isomorphisms.

A natural transformation is tied to a relation between the
argument-places of the functor in the source and the
argument-places of the functor in the target. This relation
corresponds to a relation between occurrences of letters in
formulae, and in composing natural transformations we compose
these relations. With dinatural transformations the matter is more
complicated, and composition poses particular problems (see
\cite{P03}). In this paper we deal with natural transformations.
Our general notion of coherence does not, however, presuppose
naturality and dinaturality.

Our notion of a coherence result is one that covers Mac Lane's and
Kelly's coherence results mentioned above, but it is more general.
We call coherence a result that tells us that there is a faithful
functor $G$ from a category \eS\ freely generated in a certain
class of categories to a ``manageable'' category \eM. This calls
for some explanation.

It is desirable, though perhaps not absolutely necessary, that the
functor $G$ be {\it structure-preserving}, which means that it
preserves structure at least up to isomorphism. In all coherence
results we will consider here, the functor $G$ will preserve
structure strictly, i.e.\ ``on the nose''. The categories \eS\ and
\eM\ will be in the same class of categories, and $G$ will be
obtained by extending in a unique way a map from the generators of
\eS\ into \eM.

The category \eM\ is {\it manageable} when equations of arrows,
i.e.\ commuting diagrams of arrows, are easier to consider in it
than in \eS. The best is if the commuting problem is obviously
decidable in \eM, while it was not obvious that it is such in \eS.

With our approach to coherence we are oriented towards solving the
commuting problem. This should be stressed because other authors
may give a more prominent place to other problems. We have used on
purpose the not very precise term ``manageable'' for the category
\eM\ to leave room for modifications of our notion of coherence,
which would be oriented towards solving another problem than the
commuting problem.

In this paper, the manageable category \eM\ will be the category
\emph{Rel} with arrows being relations between occurrences of
letters in formulae. In \cite{DP04} and elsewhere we have taken
\emph{Rel} to be the category of relations between finite
ordinals, which is not essentially different from what we do in
this paper. The previous category \emph{Rel} is the skeleton of
the new one. We have mentioned above the connection between
\emph{Rel} and natural transformations. The commuting problem in
\emph{Rel} is obviously decidable.

The freely generated category \eS\ will be the bicartesian
category, i.e. category with all finite products and coproducts,
freely generated by a set of objects, or a related category of
that kind. The generating set of objects may be conceived as a
discrete category. In our understanding of coherence, replacing
this discrete generating category by an arbitrary category would
prevent us to solve coherence---simply because the commuting
problem in the arbitrary generating category may be undecidable.
Far from having more general, stronger, results if the generating
category is arbitrary, we may end up by having no result at all.

The categories \eS\ in this paper are built ultimately out of {\it
syntactic} material, as logical systems are built. Categorists are
not inclined to formulate their coherence results in the way we
do---in particular, they do not deal often with syntactically
built categories. If, however, more involved and more abstract
formulations of coherence that may be found in the literature (for
early references on this matter see \cite{K72}) have practical
consequences for solving the commuting problem, our way of
formulating coherence has these consequences as well.

That there is a faithful structure-preserving functor $G$ from the
syntactical category \eS\ to the manageable category \eM\ means
that for all arrows $f$ and $g$ of \eS\ with the same source and
the same target we have
\[
f=g {\mbox{\rm{ in }}} \eS\ {\mbox {\rm{ iff }}}\; Gf=Gg
{\mbox{\rm{ in }}} \eM.
\]
The direction from left to right in this equivalence is contained
in the functoriality of $G$, while the direction from right to
left is faithfulness proper.

If \eS\ is conceived as a syntactical system, while \eM\ is a
model, the faithfulness equivalence we have just stated is like a
completeness result in logic. The left-to-right direction, i.e.\
functoriality, is soundness, while the right-to-left direction,
i.e.\ faithfulness, is completeness proper.

If $G$ happens to be one-one on objects, then we obtain that \eS\
is isomorphic to a subcategory of \eM---namely, its image under
$G$ in \eM. We will have such a situation in this paper, where $G$
will be identity on objects.

In this paper we will separate coherence results involving
terminal objects and initial objects from those not involving
them. These objects cause difficulties, and the statements and
proofs of the coherence results gain by having these difficulties
kept apart.

\section{\large\bf Coherence and proof theory}
If one envisages a deductive system as a graph whose nodes are
formulae:

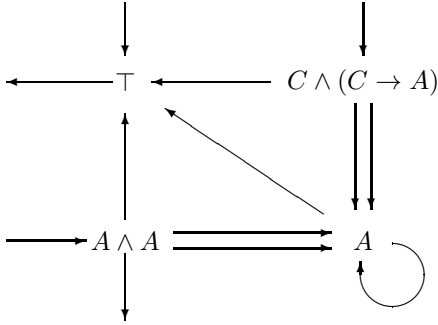
\begin{center}
\begin{picture}(150,120)
\put(45,30){\makebox(0,0){$A\wedge A$}}
\put(45,90){\makebox(0,0){$\top$}}
\put(135,30){\makebox(0,0){$A$}}
\put(135,90){\makebox(0,0){$C\wedge (C\rightarrow A)$}}

\put(0,30){\vector(1,0){30}} \put(63,27){\vector(1,0){60}}
\put(63,33){\vector(1,0){60}} \put(45,25){\vector(0,-1){25}}
\put(45,38){\vector(0,1){40}} \put(132,82){\vector(0,-1){40}}
\put(138,82){\vector(0,-1){40}} \put(120,40){\vector(-3,2){60}}
\put(40,90){\vector(-1,0){40}} \put(100,90){\vector(-1,0){45}}
\put(45,120){\vector(0,-1){20}} \put(135,120){\vector(0,-1){20}}

\put(146,17){\oval(24,24)[r]} \put(146,17){\oval(24,24)[bl]}
\put(134,17){\vector(0,1){5}}

\end{picture}
\end{center}

\noindent and whose arrows are derivations from the sources
understood as premises to the targets understood as conclusions,
then equality of derivations usually transforms this deductive
system into a category of a particular brand. This category has a
structure induced by the connectives of the deductive system.
Although equality of derivation is dictated by logical concerns,
usually the categories we end up with are of a kind that
categorists have already introduced for their own reasons. The
prime example here is given by the deductive system for the
conjunction-implication fragment of intuitionistic propositional
logic. After derivations in this deductive system are equated
according to ideas about normalization of derivations that stem
from Gentzen, one obtains the cartesian closed category $\cal K$
freely generated by a set of propositional letters (see
\cite{LS86} for the notion of cartesian closed category).

Equality of proofs in intuitionistic logic has not led up to now
to a coherence result---a coherence theorem is not forthcoming for
cartesian closed categories. If we take that the model category
$\cal M$ is a category whose arrows are graphs like the graphs of
\cite{KML71}, then we do not have a faithful functor $G$ from the
free cartesian closed category $\cal K$ to $\cal M$. We will now
explain why $G$ is not even a functor.

If $\eta_{p,q}$ is the canonical arrow from $q$ to $p\rightarrow
(p\kon q)$, where $A\rightarrow B$ and $A\kon B$ stand for $B^A$
and $A\times B$ respectively, while $\hat{w}_A$ is the diagonal
arrow from $A$ to $A\kon A$, then $G(\hat{w}_{p\rightarrow (p\kon
q)}\cirk\eta_{p,q})$:

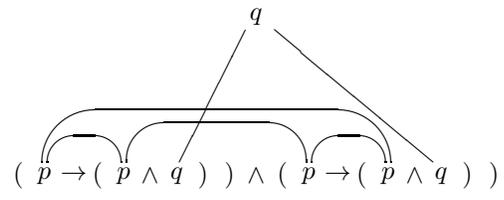
\begin{center}
\begin{picture}(180,80)
\put(10,10){\makebox(0,0){$($}} \put(20,10){\makebox(0,0){$p$}}
\put(31,10){\makebox(0,0){$\rightarrow$}}
\put(40,10){\makebox(0,0){$($}} \put(50,10){\makebox(0,0){$p$}}
\put(60,10){\makebox(0,0){$\kon$}} \put(70,10){\makebox(0,0){$q$}}
\put(80,10){\makebox(0,0){$)$}} \put(90,10){\makebox(0,0){$)$}}
\put(100,10){\makebox(0,0){$\kon$}}
\put(110,10){\makebox(0,0){$($}} \put(120,10){\makebox(0,0){$p$}}
\put(131,10){\makebox(0,0){$\rightarrow$}}
\put(140,10){\makebox(0,0){$($}} \put(150,10){\makebox(0,0){$p$}}
\put(160,10){\makebox(0,0){$\kon$}}
\put(170,10){\makebox(0,0){$q$}} \put(180,10){\makebox(0,0){$)$}}
\put(190,10){\makebox(0,0){$)$}}

\put(100,70){\makebox(0,0){$q$}}

\put(35,15){\oval(28,20)[t]}

\put(85,15){\oval(132,40)[t]}

\put(135,15){\oval(28,20)[t]}

\put(85,15){\oval(68,30)[t]}

\put(71,15){\line(1,2){25}}

\put(167,15){\line(-6,5){60}}

\end{picture}
\end{center}

\noindent which is obtained from

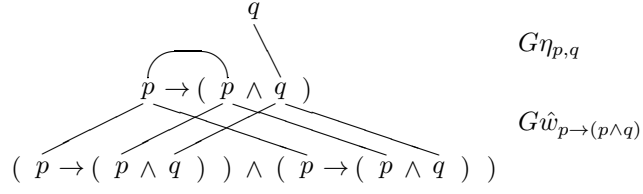
\begin{center}
\begin{picture}(180,80)
\put(10,10){\makebox(0,0){$($}} \put(20,10){\makebox(0,0){$p$}}
\put(31,10){\makebox(0,0){$\rightarrow$}}
\put(40,10){\makebox(0,0){$($}} \put(50,10){\makebox(0,0){$p$}}
\put(60,10){\makebox(0,0){$\kon$}} \put(70,10){\makebox(0,0){$q$}}
\put(80,10){\makebox(0,0){$)$}} \put(90,10){\makebox(0,0){$)$}}
\put(100,10){\makebox(0,0){$\kon$}}
\put(110,10){\makebox(0,0){$($}} \put(120,10){\makebox(0,0){$p$}}
\put(131,10){\makebox(0,0){$\rightarrow$}}
\put(140,10){\makebox(0,0){$($}} \put(150,10){\makebox(0,0){$p$}}
\put(160,10){\makebox(0,0){$\kon$}}
\put(170,10){\makebox(0,0){$q$}} \put(180,10){\makebox(0,0){$)$}}
\put(190,10){\makebox(0,0){$)$}}

\put(60,40){\makebox(0,0){$p$}}
\put(71,40){\makebox(0,0){$\rightarrow$}}
\put(80,40){\makebox(0,0){$($}} \put(90,40){\makebox(0,0){$p$}}
\put(100,40){\makebox(0,0){$\kon$}}
\put(110,40){\makebox(0,0){$q$}} \put(120,40){\makebox(0,0){$)$}}

\put(100,70){\makebox(0,0){$q$}}

\put(100,65){\line(1,-2){10}} \put(75,45){\oval(30,20)[t]}

\put(58,35){\line(-2,-1){38}} \put(88,35){\line(-2,-1){38}}
\put(108,35){\line(-2,-1){38}} \put(62,35){\line(3,-1){58}}
\put(92,35){\line(3,-1){58}} \put(112,35){\line(3,-1){58}}

\put(200,25){$G\hat{w}_{p\rightarrow (p\kon q)}$}

\put(200,55){$G\eta_{p,q}$}

\end{picture}
\end{center}

\noindent is different from $G((\eta_{p,q}\kon\eta_{p,q})\cirk
\hat{w}_q)$:

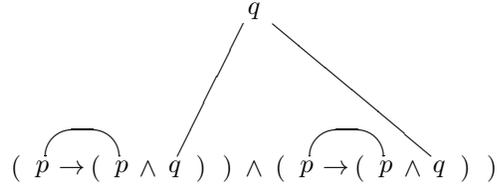
\begin{center}
\begin{picture}(180,80)
\put(10,10){\makebox(0,0){$($}} \put(20,10){\makebox(0,0){$p$}}
\put(31,10){\makebox(0,0){$\rightarrow$}}
\put(40,10){\makebox(0,0){$($}} \put(50,10){\makebox(0,0){$p$}}
\put(60,10){\makebox(0,0){$\kon$}} \put(70,10){\makebox(0,0){$q$}}
\put(80,10){\makebox(0,0){$)$}} \put(90,10){\makebox(0,0){$)$}}
\put(100,10){\makebox(0,0){$\kon$}}
\put(110,10){\makebox(0,0){$($}} \put(120,10){\makebox(0,0){$p$}}
\put(131,10){\makebox(0,0){$\rightarrow$}}
\put(140,10){\makebox(0,0){$($}} \put(150,10){\makebox(0,0){$p$}}
\put(160,10){\makebox(0,0){$\kon$}}
\put(170,10){\makebox(0,0){$q$}} \put(180,10){\makebox(0,0){$)$}}
\put(190,10){\makebox(0,0){$)$}}

\put(100,70){\makebox(0,0){$q$}}

\put(35,15){\oval(28,20)[t]}

\put(135,15){\oval(28,20)[t]}

\put(71,15){\line(1,2){25}}

\put(167,15){\line(-6,5){60}}

\end{picture}
\end{center}

\noindent which is obtained from

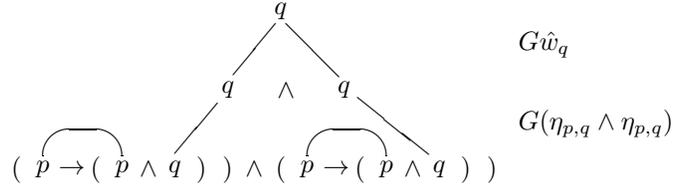
\begin{center}
\begin{picture}(400,80)(125,0)

\put(210,10){\makebox(0,0){$($}} \put(220,10){\makebox(0,0){$p$}}
\put(231,10){\makebox(0,0){$\rightarrow$}}
\put(240,10){\makebox(0,0){$($}} \put(250,10){\makebox(0,0){$p$}}
\put(260,10){\makebox(0,0){$\kon$}}
\put(270,10){\makebox(0,0){$q$}} \put(280,10){\makebox(0,0){$)$}}
\put(290,10){\makebox(0,0){$)$}}
\put(300,10){\makebox(0,0){$\kon$}}
\put(310,10){\makebox(0,0){$($}} \put(320,10){\makebox(0,0){$p$}}
\put(331,10){\makebox(0,0){$\rightarrow$}}
\put(340,10){\makebox(0,0){$($}} \put(350,10){\makebox(0,0){$p$}}
\put(360,10){\makebox(0,0){$\kon$}}
\put(370,10){\makebox(0,0){$q$}} \put(380,10){\makebox(0,0){$)$}}
\put(390,10){\makebox(0,0){$)$}}

\put(290,40){\makebox(0,0){$q$}}
\put(312,40){\makebox(0,0){$\kon$}}
\put(334,40){\makebox(0,0){$q$}}

\put(310,70){\makebox(0,0){$q$}}

\put(308,65){\line(-5,-6){16}} \put(312,65){\line(1,-1){20}}

\put(235,15){\oval(30,20)[t]} \put(335,15){\oval(30,20)[t]}

\put(286,35){\line(-5,-6){16}} \put(339,37){\line(5,-4){27}}

\put(400,25){$G(\eta_{p,q}\kon\eta_{p,q})$}

\put(400,55){$G\hat{w}_q$}

\end{picture}
\end{center}

\noindent So, if $\hat{w}$ is a natural transformation, then $G$
is not a functor. The naturality of $\hat{w}$, and other arrows of
that kind, tied to structural rules ($\hat{w}$ is tied to
contraction, and $\hat{k}^1$ below to thinning), is desirable
because it corresponds to the permuting of these rules in a
cut-elimination or normalization procedure.

Dually, if $\varepsilon_{p,q}$ is the canonical arrow from
$p\kon(p\rightarrow q)$ to $q$, and $\hat{k}^1_{A,B}$ is the first
projection from $A\kon B$ to $A$, then
$G(\hat{k}^1_{r,q}\cirk(\mj_r\kon\varepsilon_{p,q}))$:

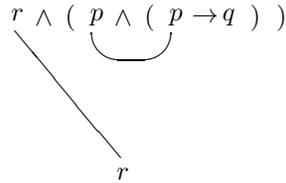
\begin{center}
\begin{picture}(120,80)

\put(50,11){\makebox(0,0){$r$}}

\put(10,71){\makebox(0,0){$r$}}

\put(20,70){\makebox(0,0){$\kon$}}

\put(30,70){\makebox(0,0){$($}}

\put(40,70){\makebox(0,0){$p$}}

\put(50,70){\makebox(0,0){$\kon$}}

\put(60,70){\makebox(0,0){$($}}

\put(70,70){\makebox(0,0){$p$}}

\put(81,70){\makebox(0,0){$\rightarrow$}}

\put(90,70){\makebox(0,0){$q$}}

\put(100,70){\makebox(0,0){$)$}}

\put(110,70){\makebox(0,0){$)$}}

\put(49,17){\line(-5,6){40}}

\put(53,64){\oval(30,20)[b]}

\end{picture}
\end{center}

\noindent which is obtained from

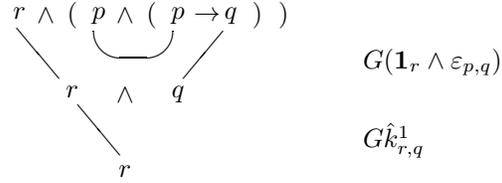
\begin{center}
\begin{picture}(120,80)

\put(50,11){\makebox(0,0){$r$}}

\put(30,41){\makebox(0,0){$r$}}

\put(50,40){\makebox(0,0){$\kon$}}

\put(70,40){\makebox(0,0){$q$}}

\put(10,71){\makebox(0,0){$r$}}

\put(20,70){\makebox(0,0){$\kon$}}

\put(30,70){\makebox(0,0){$($}}

\put(40,70){\makebox(0,0){$p$}}

\put(50,70){\makebox(0,0){$\kon$}}

\put(60,70){\makebox(0,0){$($}}

\put(70,70){\makebox(0,0){$p$}}

\put(81,70){\makebox(0,0){$\rightarrow$}}

\put(90,70){\makebox(0,0){$q$}}

\put(100,70){\makebox(0,0){$)$}}

\put(110,70){\makebox(0,0){$)$}}

\put(48,17){\line(-5,6){16}}

\put(25,46){\line(-5,6){16}}

\put(72,46){\line(5,6){15}}

\put(53,64){\oval(30,20)[b]}

\put(140,20){$G\hat{k}^1_{r,q}$}

\put(140,50){$G(\mj_r\kon\varepsilon_{p,q})$}

\end{picture}
\end{center}

\noindent is different from $G\hat{k}^1_{r,p\kon(p\rightarrow
q)}$:

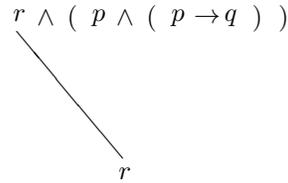
\begin{center}
\begin{picture}(120,80)

\put(50,11){\makebox(0,0){$r$}}

\put(10,71){\makebox(0,0){$r$}}

\put(20,70){\makebox(0,0){$\kon$}}

\put(30,70){\makebox(0,0){$($}}

\put(40,70){\makebox(0,0){$p$}}

\put(50,70){\makebox(0,0){$\kon$}}

\put(60,70){\makebox(0,0){$($}}

\put(70,70){\makebox(0,0){$p$}}

\put(81,70){\makebox(0,0){$\rightarrow$}}

\put(90,70){\makebox(0,0){$q$}}

\put(100,70){\makebox(0,0){$)$}}

\put(110,70){\makebox(0,0){$)$}}

\put(49,17){\line(-5,6){40}}

\end{picture}
\end{center}

\noindent So, if $\hat{k}^1$ is a natural transformation, then $G$
is not a functor. The faithfulness of $G$ fails because of a
counterexample in \cite{S75}, involving a natural number object in
\emph{Set} and the successor function. This does not exclude that
with a more sophisticated model category $\cal M$ we might still
be able to obtain coherence for cartesian closed categories (for
an attempt along these lines see \cite{PD}).

Equality of proofs in classical logic may, however, lead to
coherence with respect to model categories that catch up to a
point the idea of \emph{generality} of proofs. Such is in
particular the category \emph{Rel} mentioned in the preceding
section, whose arrows are relations between occurrences of
propositional letters in the premises and conclusions. The idea
that generality of proofs may serve as a criterion for identity of
proofs stems from Lambek's pioneering papers in categorial proof
theory of the late 1960s (see \cite{LS86} for references). This
criterion says, roughly, that two derivations represent the same
proof when their generalizations with respect to diversification
of variables (without changing the rules of inference) produce
derivations with the same source and target, up to a renaming of
variables.

Although coherence with respect to \emph{Rel} is related to
generality, it is not exactly that. The question is should
$G\hat{w}_p$ be the relation in the left one or in the right one
of the following two diagrams:

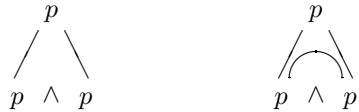
\begin{center}
\begin{picture}(140,40)

\put(120,10){\oval(20,20)[t]}

\put(6,10){\line(1,2){9}} \put(34,10){\line(-1,2){9}}
\put(106,10){\line(1,2){9}} \put(134,10){\line(-1,2){9}}

\put(7,0){\makebox(0,0)[b]{$p$}}
\put(20,1){\makebox(0,0)[b]{$\kon$}}
\put(33,0){\makebox(0,0)[b]{$p$}}

\put(107,0){\makebox(0,0)[b]{$p$}}
\put(120,1){\makebox(0,0)[b]{$\kon$}}
\put(133,0){\makebox(0,0)[b]{$p$}}

\put(20,33){\makebox(0,0)[b]{$p$}}
\put(120,33){\makebox(0,0)[b]{$p$}}

\end{picture}
\end{center}

\noindent The second option, induced by dealing with equivalence
relations, or by connecting all letters that must remain the same
in generalizing proofs (see \cite{DP03a} and \cite{DP03b}), would
lead to abolishing the naturality of $\hat{w}$. For example, in
the following instance of the naturality equation for $\hat{w}$:
\[
\hat{w}_p\cirk\check{\kappa}_p\;=(\check{\kappa}_p\kon\check{\kappa}_p)\cirk
\hat{w}_\bot
\]
for $\check{\kappa}_p$ being the unique arrow from the initial
object $\bot$ to $p$, we do not have that
$G(\hat{w}_p\cirk\check{\kappa}_p)$ is equal to
$G((\check{\kappa}_p\kon\check{\kappa}_p)\cirk \hat{w}_\bot)$:

\begin{center}
\begin{picture}(200,60)

\put(20,10){\oval(20,20)[t]}

\put(6,10){\line(1,2){9}} \put(34,10){\line(-1,2){9}}

\put(7,0){\makebox(0,0)[b]{$p$}}
\put(20,1){\makebox(0,0)[b]{$\kon$}}
\put(33,0){\makebox(0,0)[b]{$p$}}

\put(20,33){\makebox(0,0)[b]{$p$}}
\put(20,55){\makebox(0,0)[b]{$\bot$}}

\put(147,0){\makebox(0,0)[b]{$p$}}
\put(160,1){\makebox(0,0)[b]{$\kon$}}
\put(173,0){\makebox(0,0)[b]{$p$}}
\put(147,33){\makebox(0,0)[b]{$\bot$}}
\put(160,33){\makebox(0,0)[b]{$\kon$}}
\put(173,33){\makebox(0,0)[b]{$\bot$}}
\put(160,55){\makebox(0,0)[b]{$\bot$}}

\put(45,20){\makebox(0,0)[l]{$G\hat{w}_p$}}
\put(45,48){\makebox(0,0)[l]{$G\check{\kappa}_p$}}
\put(185,48){\makebox(0,0)[l]{$G\hat{w}_\bot$}}
\put(185,20){\makebox(0,0)[l]{$G(\check{\kappa}_p\kon\check{\kappa}_p)$}}

\end{picture}
\end{center}

\noindent We obtain similarly that $\kappa$ cannot be natural.

It is shown in \cite{DP04} that coherence with respect to the
model category \emph{Rel} could justify plausibly equality of
derivations in various systems of propositional logic, including
classical propositional logic. The goal of that book was to
explore the limits of coherence with respect to the model category
\emph{Rel}. This does not exclude that other coherence results may
involve other model categories, and, in particular, with a model
category different from \emph{Rel}, classical propositional logic
may induce a different notion of Boolean category than the one
introduced in Chapter 14 of \cite{DP04}. That notion of Boolean
category was not motivated \emph{a priori}, but was dictated by
coherence with respect to \emph{Rel}. The definition of that
notion was however not given via coherence, but via an equational
axiomatization. We take such definitions as being proper axiomatic
definitions.

We could easily define nonaxiomatically a notion of Boolean
category with respect to graphs of the Kelly-Mac Lane kind (see
\cite{KML71}). Equality of graphs would dictate what arrows are
equal. In this notion, conjunction would not be a product, because
the diagonal arrows and the projections would not make natural
transformations (see above), and, analogously, disjunction would
not be a coproduct (cf.\ \cite{DP04}, Section 14.3.) The resulting
notion of Boolean category would not be trivial---the freely
generated categories of that kind would not be preorders---, but
its nonaxiomatic definition would be trivial. There might exist a
nontrivial equational axiomatic definition of this notion. Finding
such a definition is an open problem.

We are looking for nontrivial axiomatic definitions because such
definitions give information about the combinatorial building
blocks of our notions, as Reidemeister moves give information
about the combinatorial building blocks of knot equivalence. Our
axiomatic equational definition of Boolean category in \cite{DP04}
is of the nontrivial, combinatorially informative, kind. Coherence
of these Boolean categories with respect to \emph{Rel} is a
theorem, whose proof in \cite{DP04} requires considerable effort.

Another analogous example is provided by the notion of monoidal
category, which was introduced in a not entirely axiomatic way,
via coherence, by B\' enabou in \cite{Ben63}, and in the axiomatic
way, such as we favour, by Mac Lane in \cite{ML63}. For B\'
enabou, coherence is built into the definition, and for Mac Lane
it is a theorem. One could analogously define the theorems of
classical propositional logic as being the tautologies (this is
done, for example, in \cite{CK73}, Sections 1.2-3), in which case
completeness would not be a theorem, but would be built into the
definition.

In this paper we prove coherence for categories that formalize
equality of proofs in classical and intuitionistic
conjunctive-disjunctive logic without distribution of conjunction
over disjunction. This fragment of logic also covers the additive
connectives of linear and other substructural logics (where
distribution anyway should not be assumed). When to this fragment
we add the true and absurd propositional constants matters become
more complicated, and we do not know how to prove unrestricted
coherence in all cases.

\section{\large\bf Lattice categories}
In the remaining sections of this paper we deal with coherence
with respect to \emph{Rel} for categories with a double cartesian
structure, i.e.\ with finite products and finite coproducts. We
take this as a categorification of the notion of lattice. As
before, we distinguish cases with and without special objects,
which are here the empty product and the empty coproduct, i.e.\
the terminal and initial objects. Categories with all finite
products and coproducts, including the empty ones, are usually
called \emph{bicartesian} categories (see \cite{LS86}). Categories
with all nonempty finite products and coproducts are called
\emph{lattice} categories in \cite{DP04}. The results presented
here are adapted from \cite{DP01}, \cite{DP02}, the revised
version of \cite{DP01a} and \cite{DP04}, Chapter~9.

We pay particular attention to questions of maximality, i.e.\ to
the impossibility of extending our axioms without collapse into
preorder, and hence triviality. This maximality is a kind of
syntactical completeness. (The sections on maximality improve upon
results reported in \cite{DP01}, \cite{DP02} and \cite{DP01a}, and
are taken over from \cite{DP04}, Chapter~9.)

Our techniques are partly based on a composition elimination for
conjunctive logic, related to normalization in natural deduction,
and on a simple composition elimination for
conjunctive-disjunctive logic, implicit in Gentzen's cut
elimination.

We define now the category $\ml$ built out of syntactic material.
The objects of the category $\ml$ are the formulae of the
propositional language $\eL$, generated out of a set of infinitely
many propositional letters, for which we use $p$, $q$, $r,\ldots$,
sometimes with indices, with the binary connectives $\kon$ and
$\dis$, for which we use $\!\ks\!$. For formulae we use $A$, $B$,
$C,\ldots$, sometimes with indices.

To define the arrows of $\ml$, we define first inductively a set
of expressions called the \emph{arrow terms} of $\ml$. Every arrow
term will have a \emph{type}, which is an ordered pair of formulae
of $\eL_{\kon}$. We write $f\!:A\vdash B$ when the arrow term $f$
is of type $(A,B)$. Here $A$ is the \emph{source}, and $B$ the
\emph{target} of $f$. For arrow terms we use $f$, $g$, $h,\ldots$,
sometimes with indices. Intuitively, the arrow term $f$ is the
code of a derivation of the conclusion $B$ from the premise $A$
(which explains why we write $\vdash$ instead of $\str$).

For all formulae $A$, $B$ and $C$ of $\eL$ the following
\emph{primitive arrow terms}:
\begin{tabbing}
\hspace{13em}$\mj_A\!: A\vdash A$,
\\[1ex]
\mbox{\hspace{6em}}\=$\hat{w}_{A}\!:A\vdash A\kon
A$,\hspace{4em}\= $\check{w}_{A}\!:A\dis A \vdash A$,
\\[1ex]
\>$\hat{k}^i_{A_1,A_2}\!:A_1\kon A_2\vdash
A_i$,\>$\check{k}^i_{A_1,A_2}\!:A_i\vdash A_1\dis A_2$,
\hspace{1em} for $i\in\{1,2\}$,
\end{tabbing}
are arrow terms. (Intuitively, these are the axioms of our logic
with the codes of their trivial derivations.)

Next we have the following inductive clauses:
\begin{itemize}
\item[]if ${f\!:A\vdash B}$ and ${g\!:B\vdash C}$ are arrow
terms,\\ then ${(g\cirk f)\!:A\vdash C}$ is an arrow
term;\vspace{-1ex} \item[]if ${f_1\!:A_1\vdash B_1}$ and
${f_2\!:A_2\vdash B_2}$ are arrow terms,\\ then ${(f_1\ks
f_2)\!:A_1\ks A_2\vdash B_1\ks B_2}$ is an arrow term.
\end{itemize}
(Intuitively, the operations on arrow terms $\cirk$ and $\!\ks\!$
are codes of the rules of inference of our logic.) This defines
the arrow terms of $\ml$. As we do usually with formulae, we will
omit the outermost parentheses of arrow terms.

We stipulate first that all the instances of ${f=f}$ and of the
following equations are equations of $\ml$:
\begin{tabbing}
\hspace{1em}\emph{categorial equations}:\\*[1ex]
\mbox{\hspace{2em}}\= $({\mbox{{\it cat}~1}})$\hspace{2em}\=
$f\cirk \mj_A=\mj_B\cirk f=f\!:A\vdash B$,
\\*[1ex]
\> $({\mbox{{\it cat}~2}})$\> $h\cirk (g\cirk f)=(h\cirk g)\cirk
f$,
\\[1.5ex]
\hspace{1em}\emph{bifunctorial equations}:\\*[1ex]
  \> $(\!\ks\, 1)$\> $\mj_A\ks\mj_B=\mj_{A\kst B}$,
\\*[1ex]
\> $(\!\ks\, 2)$\> $(g_1\cirk f_1)\ks(g_2\cirk f_2)=(g_1\ks
g_2)\cirk(f_1\ks f_2)$,
\\[2ex]
\hspace{1em}\emph{naturality equations}:\hspace{1em}for
$f\!:A\vdash B$ and $f_i\!:A_i\vdash B_i$, where $i\in\{1,2\}$,
\\*[1ex]
\> $\mbox{($\hat{w}$ {\it nat})}$\>  $(f\kon
f)\cirk\hat{w}_{A}\:=\:\hat{w}_{B}\cirk f$,\\*[1ex]
\>$\mbox{($\check{w}$ {\it nat})}$\> $f \cirk\check{w}_{A}\:=\:
\check{w}_{B}\cirk (f\dis f)$,
\\[1ex]
\> $\mbox{($\hat{k}^i$ {\it nat})}$\>
$f_i\cirk\hat{k}^i_{A_1,A_2}\:=\:\hat{k}^i_{B_1,B_2}\cirk (f_1\kon
f_2)$,\\*[1ex] \> $\mbox{($\check{k}^i$ {\it nat})}$\> $(f_1\dis
f_2)\cirk\check{k}^i_{A_1,A_2}\:=\:\check{k}^i_{B_1,B_2}\cirk
f_i$,
\\[1.5ex]
\hspace{1em}\emph{triangular equations}:\hspace{1em}for
$i\in\{1,2\}$,\\*[1ex]

\> ($\hat{w}\hat{k}$) \>
$\hat{k}^{i}_{A,A}\cirk\hat{w}_A\,$\=$=\mj_A$,\\*[1ex] \>
($\check{w}\check{k}$) \>
$\check{w}_A\cirk\check{k}^{i}_{A,A}$\>$=\mj_A$,\\[1ex]

\> ($\hat{w}\hat{k}\hat{k}$)\> $(\hat{k}^1_{A,B}\kon
\hat{k}^2_{A,B})\cirk \hat{w}_{A\kon B}\,$\=$= \mj_{A\kon
B}$,\\*[1ex]

\> ($\check{w}\check{k}\check{k}$)\> $\check{w}_{A\dis
B}\cirk(\check{k}^1_{A,B}\dis \check{k}^2_{A,B})$\>$= \mj_{A\dis
B}$.
\end{tabbing}

This concludes the list of axiomatic equations stipulated for \ml.
To define all the equations of \ml\ it remains only to say that
the set of these equations is closed under symmetry and
transitivity of equality and under the rules
\[
(\!\cirk\;\mbox{\it cong})\quad \f{f=f'\quad \quad \quad g=g'}
{g\cirk f=g'\cirk f'}\hspace{5em}(\!\ks\;\mbox{\it cong})\quad
\f{f_1=f_1'\quad \quad \quad f_2=f_2'} {f_1\ks f_2=f_1'\ks f_2'}
\]

On the arrow terms of \ml\ we impose the equations of \ml. This
means that an arrow of \ml\ is an equivalence class of arrow terms
of \ml\ defined with respect to the smallest equivalence relation
such that the equations of \ml\ are satisfied (see \cite{DP04},
Section 2.3, for details).

The kind of category for which \ml\ is the one freely generated
out of the set of propositional letters (which may be understood
as a discrete category) we call \emph{lattice} category (see
\cite{DP04}, Section 9.4, for a precise definition). Usually, such
categories would be called categories with finite nonempty
products and coproducts. The objects of a lattice category that is
a partial order make a lattice.

\section{\large\bf The functor $G$}
The objects of the category \emph{Rel} are the objects of \ml,
i.e.\ the formulae of $\eL$. An arrow $R\!:A\vdash B$ of
\emph{Rel} is a set of ordered pairs $(x,y)$ such that $x$ is an
occurrence of a propositional letter in the formula $A$ and $y$ is
an occurrence of a propositional letter in the formula $B$; in
other words, arrows are binary relations between the sets of
occurrences of propositional letters in formulae. We write either
$(x,y)\in R$ or $xRy$, as usual. In this category,
$\mj_A\!:A\vdash A$ is the identity relation, i.e.\ identity
function, that assigns to every occurrence of a propositional
letter in $A$ that same occurrence. In $\eL$ there are no formulae
in which no propositional letter occurs, but where we have such
formulae (as in the language $\eL_{\top,\bot}$ considered later in
this paper), the empty set of ordered pairs corresponds to
$\mj_A\!:A\vdash A$ if no propositional letter occurs in $A$. The
empty relation is the identity relation on the empty set.

For $R_1\!:A\vdash B$ and $R_2\!:B\vdash C$, the set of ordered
pairs of the composition $R_2\cirk R_1\!:A\vdash C$ is
$\{(x,y)\mid \exists z(xR_1 z\;\mbox{\it and}\;zR_2 y)\}$. Let
$x_j(A)$ be the $j$-th occurrence of a propositional letter in $A$
counting from the left, and let $|A|$ be the number of occurrences
of propositional letters in $A$ (so $1\leq j\leq |A|$). For
$R_i\!:A_i\vdash B_i$, with $i\in\{1,2\}$, the set of ordered
pairs of $R_1 \ks R_2\!:A_1\ks A_2\vdash B_1\ks B_2$, for
$\!\ks\!\in\{\kon,\dis\}$, is the disjoint union of the following
two sets:
\begin{tabbing}
\hspace{5em}\=$\{(x_j(A_1\ks A_2),x_k(B_1\ks B_2))\mid
(x_j(A_1),x_k(B_1))\in R_1\}$,\\*[1ex]

\>$\{(x_{j+|A_1|}(A_1\ks A_2),x_{k+|B_1|}(B_1\ks B_2))\mid
(x_j(A_2),x_k(B_2))\in R_2\}$.
\end{tabbing}
With the operation on objects that corresponds to the binary
connective $\!\ks\!$, this operation $\!\ks\!$ on arrows gives a
biendofunctor in \emph{Rel}.

In \emph{Rel} we have the relations $G\hat{w}_{A}\!:A\vdash A\kon
A$, $G\check{w}_{A}\!:A\dis A \vdash A$,
$G\hat{k}^i_{A_1,A_2}\!:A_1\kon A_2\vdash A_i$, and
$G\check{k}^i_{A_1,A_2}\!:A_i\vdash A_1\dis A_2$, for
$i\in\{1,2\}$, whose sets of ordered pairs are defined as follows:
\begin{tabbing}
\hspace{0em}\=$(x_j(A),x_k(A\kon A))\in G\hat{w}_A$ iff
$(x_k(A\dis A),x_j(A))\in G\check{w}_A$ iff $j\equiv k$ (mod
$|A|$);
\\[1ex]
\>$(x_j(A_1\kon A_2),x_k(A_1))\in G\hat{k}^{1}_{A_1,A_2}$ iff
$(x_k(A_1),x_j(A_1\dis A_2))\in G\check{k}^{1}_{A_1,A_2}$ iff
$j=k$;
\\[1ex]
\>$(x_j(A_1\kon A_2),x_k(A_2))\in G\hat{k}^{2}_{A_1,A_2}$ iff
$(x_k(A_2),x_j(A_1\dis A_2))\in G\check{k}^{2}_{A_1,A_2}$
iff\\*[.5ex]\` $j=k\pl |A_1|$.
\end{tabbing}

It is not difficult to check that all these arrows of \emph{Rel}
give rise to natural transformations. This is clear from the
graphical representation of relations in \emph{Rel}. Here are a
few examples of such graphical representations, with sources
written at the top and targets at the bottom:

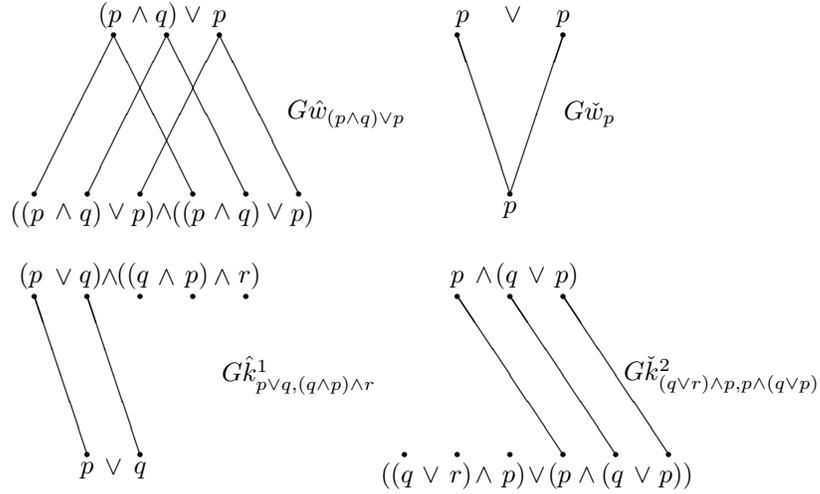
\begin{center}
\begin{picture}(260,80)
\put(10,10){\line(1,2){30}} \put(30,10){\line(1,2){30}}
\put(50,10){\line(1,2){30}} \put(70,10){\line(-1,2){30}}
\put(90,10){\line(-1,2){30}} \put(110,10){\line(-1,2){30}}

\put(190,10){\line(-1,3){20}} \put(190,10){\line(1,3){20}}

\put(10,10){\circle*{2}} \put(30,10){\circle*{2}}
\put(50,10){\circle*{2}} \put(70,10){\circle*{2}}
\put(90,10){\circle*{2}} \put(110,10){\circle*{2}}
\put(190,10){\circle*{2}}

\put(40,70){\circle*{2}} \put(60,70){\circle*{2}}
\put(80,70){\circle*{2}} \put(170,70){\circle*{2}}
\put(210,70){\circle*{2}}

\put(10,7){\makebox(0,0)[t]{$((p\;\;$}}
\put(21,6){\makebox(0,0)[t]{$\kon$}}
\put(30,7){\makebox(0,0)[t]{$\;q)$}}
\put(41,6){\makebox(0,0)[t]{$\dis$}}
\put(50,7){\makebox(0,0)[t]{$\;p)$}}
\put(59,6){\makebox(0,0)[t]{$\kon$}}
\put(71,7){\makebox(0,0)[t]{$((p\;\;$}}
\put(81,6){\makebox(0,0)[t]{$\kon$}}
\put(90,7){\makebox(0,0)[t]{$\;q)$}}
\put(101,6){\makebox(0,0)[t]{$\dis$}}
\put(110,7){\makebox(0,0)[t]{$\;p)$}}
\put(190,7){\makebox(0,0)[t]{$p$}}

\put(40,75){\makebox(0,0)[b]{$(p\;$}}
\put(50,75){\makebox(0,0)[b]{$\kon$}}
\put(60,75){\makebox(0,0)[b]{$q)$}}
\put(70,75){\makebox(0,0)[b]{$\dis$}}
\put(80,75){\makebox(0,0)[b]{$p$}}
\put(172,75){\makebox(0,0)[b]{$p$}}
\put(191,75){\makebox(0,0)[b]{$\dis$}}
\put(210,75){\makebox(0,0)[b]{$p$}}

\put(128,40){\makebox(0,0){$G\hat{w}_{(p\kon q)\dis p}$}}
\put(220,40){\makebox(0,0){$G\check{w}_p$}}

\end{picture}
\end{center}

\vspace{1ex}

\begin{center}
\begin{picture}(260,80)
\put(30,10){\line(-1,3){20}}

\put(50,10){\line(-1,3){20}}

\put(210,10){\line(-2,3){40}} \put(230,10){\line(-2,3){40}}
\put(250,10){\line(-2,3){40}}

\put(30,10){\circle*{2}}

\put(50,10){\circle*{2}}

\put(150,10){\circle*{2}}

\put(170,10){\circle*{2}}

\put(210,10){\circle*{2}}

\put(190,10){\circle*{2}}

\put(230,10){\circle*{2}}

\put(250,10){\circle*{2}}

\put(10,70){\circle*{2}} \put(30,70){\circle*{2}}
\put(50,70){\circle*{2}} \put(70,70){\circle*{2}}
\put(90,70){\circle*{2}} \put(170,70){\circle*{2}}
\put(190,70){\circle*{2}} \put(210,70){\circle*{2}}

\put(30,7){\makebox(0,0)[t]{$p$}}
\put(40,8){\makebox(0,0)[t]{$\dis$}}
\put(50,7){\makebox(0,0)[t]{$q$}}
\put(150,7){\makebox(0,0)[t]{$((q\;\;$}}
\put(160,5){\makebox(0,0)[t]{$\dis$}}
\put(170,7){\makebox(0,0)[t]{$\;r)$}}
\put(180,5){\makebox(0,0)[t]{$\kon$}}
\put(190,7){\makebox(0,0)[t]{$\;p)$}}
\put(200,5){\makebox(0,0)[t]{$\dis$}}
\put(210,7){\makebox(0,0)[t]{$(p\;$}}
\put(219,5){\makebox(0,0)[t]{$\kon$}}
\put(230,7){\makebox(0,0)[t]{$(q\;$}}
\put(240,5){\makebox(0,0)[t]{$\dis$}}
\put(250,7){\makebox(0,0)[t]{$\;\;p))$}}

\put(10,75){\makebox(0,0)[b]{$(p\;$}}
\put(21,80){\makebox(0,0)[t]{$\dis$}}
\put(30,75){\makebox(0,0)[b]{$\;q)$}}
\put(38.5,80){\makebox(0,0)[t]{$\kon$}}
\put(50,75){\makebox(0,0)[b]{$((q\;\;$}}
\put(60,80){\makebox(0,0)[t]{$\kon$}}
\put(70,75){\makebox(0,0)[b]{$\;p)$}}
\put(81,80){\makebox(0,0)[t]{$\kon$}}
\put(90,75){\makebox(0,0)[b]{$\;r)$}}
\put(170,75){\makebox(0,0)[b]{$p$}}
\put(180,75){\makebox(0,0)[b]{$\kon$}}
\put(190,75){\makebox(0,0)[b]{$(q\;$}}
\put(200,75){\makebox(0,0)[b]{$\dis$}}
\put(210,75){\makebox(0,0)[b]{$\;p)$}}

\put(110,40){\makebox(0,0){$G\hat{k}^{1}_{p\dis q,(q\kon p)\kon
r}$}} \put(270,40){\makebox(0,0){$G\check{k}^{2}_{(q\dis r) \kon
p,p\kon(q\dis p)}$}}

\end{picture}
\end{center}

For $R\!:A\vdash B$, the naturality equation
\[
(R\kon R)\cirk G\hat{w}_A=\,G\hat{w}_B\cirk R,
\]
which corresponds to the equation $\mbox{($\check{w}$ {\it nat})}$
of the preceding section, and which we take as an example, is
justified in the following manner via graphs:

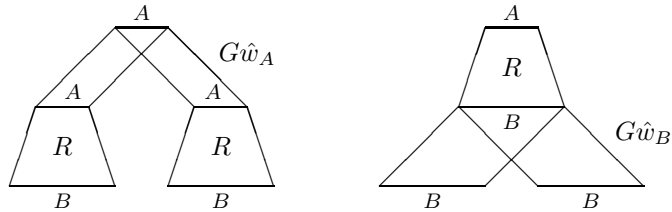
\begin{center}
\begin{picture}(260,80)
\put(10,10){\line(1,3){10}} \put(10,10){\line(1,0){40}}
\put(50,10){\line(-1,3){10}} \put(20,40){\line(1,0){20}}
\put(20,40){\line(1,1){30}} \put(40,40){\line(1,1){30}}
\put(70,10){\line(1,3){10}} \put(70,10){\line(1,0){40}}
\put(110,10){\line(-1,3){10}} \put(80,40){\line(1,0){20}}
\put(80,40){\line(-1,1){30}} \put(100,40){\line(-1,1){30}}
\put(90,50){\line(-1,1){20}} \put(50,70){\line(1,0){20}}

\put(150,10){\line(1,0){40}} \put(150,10){\line(1,1){30}}
\put(190,10){\line(1,1){30}} \put(180,40){\line(1,0){40}}
\put(210,10){\line(1,0){40}} \put(210,10){\line(-1,1){30}}
\put(250,10){\line(-1,1){30}} \put(180,40){\line(1,3){10}}
\put(220,40){\line(-1,3){10}} \put(190,70){\line(1,0){20}}

\put(30,7){\makebox(0,0)[t]{\footnotesize $B$}}
\put(60,73){\makebox(0,0)[b]{\footnotesize $A$}}
\put(90,7){\makebox(0,0)[t]{\footnotesize $B$}}
\put(30,25){\makebox(0,0){$R$}} \put(90,25){\makebox(0,0){$R$}}
\put(200,55){\makebox(0,0){$R$}}
\put(170,7){\makebox(0,0)[t]{\footnotesize $B$}}
\put(200,73){\makebox(0,0)[b]{\footnotesize $A$}}
\put(230,7){\makebox(0,0)[t]{\footnotesize $B$}}

\put(87,43){\makebox(0,0)[b]{\footnotesize $A$}}
\put(34,43){\makebox(0,0)[b]{\footnotesize $A$}}
\put(200,37){\makebox(0,0)[t]{\footnotesize $B$}}

\put(100,60){\makebox(0,0){$G\hat{w}_A$}}
\put(250,30){\makebox(0,0){$G\hat{w}_B$}}

\end{picture}
\end{center}

We can now define a functor $G$ from the category \ml\ to the
category \emph{Rel}. On objects we have that $GA$ is $A$. We have
defined $G$ above on the primitive arrow terms of \ml, and we have
\[
\begin{array}{l}
G(f\ks g)=Gf\ks Gg,
\\[.5ex]
G(g\cirk f)=Gg\cirk Gf.
\end{array}
\]
To ascertain that this defines a functor from \ml\ to \emph{Rel},
it remains to check that if $f=g$ in \ml, then $Gf=Gg$ in
\emph{Rel}, which we do by induction on the length of the
derivation of $f=g$ in \ml.

It is easy to check by induction that if for $f\!:A\vdash B$ we
have $(x_j(A),x_k(B))\in Gf$, then $x_j(A)$ and $x_k(B)$ are
occurrences of the same propositional letter.

Our first task in this paper is to show that the functor $G$ from
\ml\ to \emph{Rel} is faithful. We call this result \emph{Lattice
Coherence}, and we say that \ml\ is \emph{coherent}. Since $G$ is
identity on objects, this means that \ml\ is isomorphic to a
subcategory of \emph{Rel}.

It is clear that if \ml\ is coherent in the sense just specified,
then it is decidable whether arrow terms of \ml\ are equal in \ml.
In logical terms, one would say that the coherence of \ml\ implies
the decidability of the equational system used to define \ml. This
is because equality of arrows is clearly decidable in \emph{Rel}.
So coherence here implies a solution to the commuting problem.

\section{\large\bf Coherence for lattice categories}
We define by induction a set of terms for the arrows of \ml\ that
we call \emph{Gentzen terms}. The identity arrow terms $\mj_A$ are
Gentzen terms, and we assume that Gentzen terms are closed under
the following operations on arrow terms, besides the operation
$\cirk$, where $=_{dn}$ is read ``denotes'':
\[
\f{f_1\!:C\vdash A_1 \quad \quad \quad f_2\!:C\vdash A_2} {\langle
f_1,f_2\rangle =_{dn} (f_1\kon f_2)\cirk \hat{w}_C\!:C\vdash
A_1\kon A_2}
\]
\[
\f{g_i\!:A_i\vdash C}{\hat{K}^{i}_{A_{3-i}}g_i=_{dn}g_i\cirk
\hat{k}^i_{A_1,A_2}\!: A_1\kon A_2\vdash C}
\]
\[
\f{g_1\!:A_1\vdash C \quad \quad \quad g_2\!:A_2\vdash C}
{[g_1,g_2]=_{dn}\check{w}_C\cirk(g_1\dis g_2)\!:A_1\dis A_2\vdash
C}
\]
\[
\f{f_i\!:C\vdash A_i}
{\check{K}^{i}_{A_{3-i}}\!f_i=_{dn}\check{k}^{i}_{A_1,A_2}\cirk
f_i\!: C\vdash A_1\dis A_2}
\]

It is easy to verify that the following equations hold for Gentzen
terms (these equations can serve for an alternative formulation of
\ml):
\begin{tabbing}
\hspace{0em}\=\mbox{($\hat{K}^{}$1)}\hspace{2em}\=$g\cirk\hat{K}^{i}_A
f =\hat{K}^{i}_A (g\cirk
f)$,\hspace{5em}\=\mbox{($\check{K}^{}$1)}\hspace{2em}\=$\check{K}^{i}_A
g\cirk f =\check{K}^{i}_A (g\cirk f)$,
\\[1ex]
\>\mbox{($\hat{K}^{}$2)}\>$\hat{K}^{i}_A g\cirk\langle
f_1,f_2\rangle = g\cirk f_i$,\>\mbox{($\check{K}^{}$2)}\>
$[g_1,g_2]\cirk\check{K}^{i}_A f = g_i\cirk f$,
\\[1ex]
\>\mbox{($\hat{K}^{}$3)}\>$\langle g_1,g_2\rangle\cirk f = \langle
g_1\cirk f,g_2\cirk f\rangle$,\>\mbox{($\check{K}^{}$3)}\>$g\cirk
[f_1,f_2] = [g\cirk f_1,g\cirk f_2]$,
\\[1ex]
\>\mbox{($\hat{K}^{}$4)}\>$\mj_{A\kon B}= \langle \hat{K}^{1}_B
\mj_A,\hat{K}^{2}_A \mj_B\rangle$,\>\mbox{($\check{K}^{}$4)}\>$
\mj_{A\dis B}=
[\check{K}^{1}_B\mj_A,\check{K}^{2}_A\mj_B]$,\\[1ex]
\>\mbox{($\hat{K}^{}$5)}\>$\hat{K}^{i}_D \langle
f_1,f_2\rangle=\langle \hat{K}^{i}_D f_1,\hat{K}^{i}_D
f_2\rangle$,\>\mbox{($\check{K}^{}$5)}\>$\check{K}^{i}_D[
g_1,g_2]=[\check{K}^{i}_D g_1,\check{K}^{i}_D g_2]$,
\\[2ex]
\hspace{9em}$(\hat{K}\check{K})$\hspace{2em}$\hat{K}^{i}_C\check{K}^{j}_D
h = \check{K}^{j}_D\hat{K}^{i}_C h$,
\end{tabbing}
with appropriate types assigned to $f$, $g$, $f_i$ and $g_i$.

It is very easy to show that for every arrow term of \ml\ there is
a Gentzen term denoting the same arrow. We can prove the following
theorem for \ml.

\prop{Composition Elimination}{For every arrow term $h$ there is a
composition-free Gentzen term $h'$ such that $h=h'$.}

\dkz We find first a Gentzen term denoting the same arrow as $h$.
Take a subterm $g\cirk f$ of this Gentzen term such that both $f$
and $g$ are composition-free. We call such a subterm a {\it
topmost cut}. We show that $g\cirk f$ is equal either to a
composition-free Gentzen term or to a Gentzen term all of whose
compositions occur in topmost cuts of strictly smaller length than
the length of $g\cirk f$. The possibility of eliminating
composition in topmost cuts, and hence every composition, follows
by induction on the length of topmost cuts.

The cases where $f$ or $g$ are $\mj_A$ are taken care of by
({\mbox{{\it cat} 1}}); the cases where $f$ is $\hat{K}^{i}_A f'$
are taken care of by \mbox{($\hat{K}^{}$1)}; and the case where
$g$ is $\langle g_1,g_2\rangle$ is taken care of by
\mbox{($\hat{K}^{}$3)}.

We have next cases dual to the last two, where $g$ is
$\check{K}^{i}_Ag'$, which is taken care of by
\mbox{($\check{K}^{}$1)}, and where $f$ is $[f_1,f_2]$, which is
taken care of by \mbox{($\check{K}^{}$3)}. In the remaining cases,
if $f$ is $\langle f_1,f_2\rangle$, then $g$ is either of a form
already covered by cases above, or $g$ is $\hat{K}^{i}_A g'$, and
we apply \mbox{($\hat{K}^{}$2)}. Finally, if $f$ is
$\check{K}^{i}_A f'$, then $g$ is either of a form already covered
by cases above, or $g$ is $[g_1,g_2]$, and we apply
\mbox{($\check{K}^{}$2)}. \qed

\vspace{2ex}

\noindent Note that we use only the equations
\mbox{($\hat{K}^{}$1)}-\mbox{($\hat{K}^{}$3)} and
\mbox{($\check{K}^{}$1)}-\mbox{($\check{K}^{}$3)} in this proof
(which is taken over from \cite{DP02}, Section~3). We can then
prove the following lemma for \ml.

\prop{Invertibility Lemma for $\kon$}{Let $f\!:A_1\kon A_2\vdash
B$ be a Gentzen term. If for every $(x,y)\in Gf$ we have that $x$
is in $A_1$, then $f$ is equal to a Gentzen term of the form
${\hat{K}^{1}_{A_2}f'}$, and if for every $(x,y)\in Gf$ we have
that $x$ is in $A_2$, then $f$ is equal to a Gentzen term of the
form ${\hat{K}^{2}_{A_1}f'}$.}

\dkz By Composition Elimination for \ml, we can assume that $f$ is
com\-po\-si\-tion-free, and then we proceed by induction on the
length of the target $B$ (or on the length of $f$). If $B$ is a
letter, then $f$ must be equal in \ml\ to an arrow term of the
form ${\hat{K}^{i}_{A_{3-i}}f'}$. The condition on $Gf$ dictates
whether $i$ here is $1$ or $2$.

If $B$ is $B_1\kon B_2$ and $f$ is not of the form
${\hat{K}^{i}_{A_{3-i}}f'}$, then $f$ must be of the form $\langle
f_1,f_2\rangle$ (the condition on $Gf$ precludes that $f$ be an
identity arrow term). We apply the induction hypothesis to
${f_1\!:A_1\kon A_2\vdash B_1}$ and ${f_2\!:A_1\kon A_2\vdash
B_2}$, and use the equation \mbox{($\hat{K}^{}$5)}.

If $B$ is $B_1\dis B_2$ and $f$ is not of the form
${\hat{K}^{i}_{A_{3-i}}f'}$, then $f$ must be of the form
${\check{K}^{j}_{B_{3-j}}g}$, for ${j\in\{1,2\}}$. We apply the
induction hypothesis to ${g\!:A_1\kon A_2\vdash B_i}$, and use the
following instance of the equation $(\hat{K}\check{K})$:
\[
{\check{K}^{j}_{B_{3-j}}\hat{K}^{i}_{A_{3-i}}g'=
\hat{K}^{i}_{A_{3-i}}\check{K}^{j}_{B_{3-j}}g'}.
\]\qed

\vspace{2ex}

We have a dual Invertibility Lemma for $\dis$. We can then prove
the following result of \cite{DP02} (Section~4).

\prop{Lattice Coherence}{The functor $G$ from \ml\ to \emph{Rel}
is faithful.}

\dkz Suppose $f,g\!:A\vdash B$ are arrow terms of \ml\ and
$Gf=Gg$. We proceed by induction on the sum of the lengths of $A$
and $B$ to show that $f=g$ in \ml. If $A$ and $B$ are both
letters, then we conclude by Composition Elimination for \ml\ that
an arrow term of \ml\ of the type $A\vdash B$ exists iff $A$ and
$B$ are the same letter $p$, and we must have $f=g=\mj_p$ in \ml.
Note that we do not need here the assumption $Gf=Gg$.

If $B$ is $B_1\kon B_2$, then for $i\in\{1,2\}$ we have that
$\hat{k}^{i}_{B_1,B_2}\cirk f$ and $\hat{k}^{i}_{B_1,B_2}\cirk g$
are of type $A\vdash B_i$. We also have
\[
G(\hat{k}^{i}_{B_1,B_2}\cirk f) =G\hat{k}^{i}_{B_1,B_2}\cirk Gf
=G\hat{k}^{i}_{B_1,B_2}\cirk Gg =G(\hat{k}^{i}_{B_1,B_2}\cirk g),
\]
whence, by the induction hypothesis, we have
$\hat{k}^{i}_{B_1,B_2}\cirk f=\hat{k}^{i}_{B_1,B_2}\cirk g$ in
\ml. Then we infer
\[
\langle \hat{k}^{1}_{B_1,B_2}\cirk f,\hat{k}^{2}_{B_1,B_2}\cirk
f\rangle = \langle \hat{k}^{1}_{B_1,B_2}\cirk
g,\hat{k}^{2}_{B_1,B_2}\cirk g\rangle,
\]
from which $f=g$ follows with the help of the equations
\mbox{($\hat{K}^{}$3)} and \mbox{($\hat{K}^{}$4)}. We proceed
analogously if $A$ is $A_1\dis A_2$.

Suppose now that $A$ is $A_1\kon A_2$ or a letter, and $B$ is
$B_1\dis B_2$ or a letter, but $A$ and $B$ are not both letters.
Then by Composition Elimination for \ml\ we have that $f$ is equal
in \ml\ to an arrow term of \ml\ that is either of the form
$f'\cirk\hat{k}^{i}_{A_1,A_2}$ or of the form
$\check{k}^{i}_{B_1,B_2}\cirk f'$. Suppose
$f=f'\cirk\hat{k}^{1}_{A_1,A_2}$. Then for every $(x,y)\in Gf$ we
have $x\in GA_1$.

By the Invertibility Lemma for $\kon$, it follows that $g$ is
equal in \ml\ to an arrow term of the form
${g'\cirk\hat{k}^{1}_{A_1,A_2}}$. From ${Gf=Gg}$ we can infer
easily that ${Gf'=Gg'}$, and so by the induction hypothesis
$f'=g'$, and hence $f=g$.

We reason analogously when $f=f'\cirk\hat{k}^{2}_{A_1,A_2}$. If
$f=\:\check{k}^{i}_{B_1,B_2}\cirk f'$, then again we reason
analogously, applying the Invertibility Lemma for $\dis$. \qed

\vspace{2ex}

This proof of Lattice Coherence is simpler than a proof given in
\cite{DP02}. In the course of that previous proof one has also
coherence results for two auxiliary categories related to \ml. We
will need these categories later, but we do not need these
coherence results. For the sake of completeness, however, we
record them here too.

Let $\mlk_{\dis}$ be the category defined as \ml\ with the
difference that the primitive arrow terms $\check{w}$ and
$\check{k}^i$ are excluded, as well as the equations involving
them. The Gentzen formulation of $\mlk_{\dis}$ is obtained by
taking the operation $\dis$ on arrow terms instead of the
operations $[\;,\;]$ and $\check{K}^i$.

The category $\mld_{\kon}$ is isomorphic to
$\mlk^{\raisebox{-5pt}{\scriptsize{\it op}}}_{\dis}$. In
$\mld_{\kon}$, the $\kon$ and $\dis$ of $\mlk_{\dis}$ are
interchanged.

One can easily prove Composition Elimination for $\mlk_{\dis}$
(and hence also for $\mld_{\kon}$) by abbreviating the proof of
Composition Elimination for \ml\ above. For $\mlk_{\dis}$ we do
not have the cases where $f$ is $[f_1,f_2]$ or
$\check{K}^{i}_Af'$, but $f$ can be $f_1\dis f_2$. Then, if $g$ is
not of a form already covered by the proof above, it must be
$g_1\dis g_2$, and we apply the bifunctorial equation $(\dis\,2)$.

A composition-free arrow term of $\mlk_{\dis}$ may be reduced to a
unique normal form, which can then be used to demonstrate
coherence for $\mlk_{\dis}$, i.e.\ the fact that the functor $G$
from $\mlk_{\dis}$ to \emph{Rel} is faithful (see \cite{DP02},
Section~4).

\section{\large\bf Coherence for sesquicartesian categories}
We define now the category $\ml_{\top,\bot}$, whose definition
extends the definition of \ml\ with the terminal object $\top$ and
the initial object $\bot$, i.e.\ nullary product and coproduct.
The objects of this category are the formulae of the propositional
language $\eL_{\top,\bot}$, generated out of a set of infinitely
many propositional letters with the binary connectives $\kon$ and
$\dis$, and the nullary connectives, i.e.\ propositional
constants, $\top$ and $\bot$.

The arrow terms of $\ml_{\top,\bot}$ are defined as the arrow
terms of \ml\ save that for every object $A$ we have the
additional primitive arrow terms
\[
\hat{\kappa}_A\!:A\vdash\top, \quad\quad\quad
\check{\kappa}_A\!:\bot\vdash A,
\]
and for all arrow terms $f\!:A\vdash\top$ and $g\!:\bot\vdash A$
we have the additional axiomatic equations
\begin{tabbing}
\hspace{7em}\=\mbox{($\hat{\kappa}$)}\hspace{2em}\=$\hat{\kappa}_A=
f$,\hspace{5em}\=\mbox{($\check{\kappa}$)}\hspace{2em}\=$\check{\kappa}_A=
g$,\\[1ex]

\>$(\hat{k}^{}\bot)$\>$\hat{k}^{1}_{\bot,\bot}=\hat{k}^{2}_{\bot,\bot}$,\>
$(\check{k}^{}\top)$\>$\check{k}^{1}_{\top,\top}=\check{k}^{2}_{\top,\top}$.
\end{tabbing}

It is easy to see that with the help of the last two equations we
obtain that the pairs
\[
\begin{array}{lll}
\hat{k}^{1}_{\bot,\bot}\:=\;\hat{k}^{2}_{\bot,\bot}:\bot\kon\bot\vdash\bot
& {\mbox{\rm and}} &
\check{\kappa}_{\bot\kon\bot}\:=\;\hat{w}_{\bot}:\bot\vdash\bot\kon\bot,
\\[.5ex]
\check{k}^{1}_{\top,\top}\:=\;\check{k}^{2}_{\top,\top}:\top\vdash\top\dis\top
& {\mbox{\rm and}} &
\hat{\kappa}_{\top\dis\top}\:=\;\check{w}_{\top}:\top\dis\top\vdash\top
\end{array}
\]
are inverses of each other. This shows that every letterless
formula of $\eL_{\top,\bot}$ is isomorphic in $\ml_{\top,\bot}$
either to $\top$ or to $\bot$.

The kind of category for which $\ml_{\top,\bot}$ is the one freely
generated out of the set of propositional letters we call
\emph{dicartesian} category. The objects of a dicartesian category
that is a partial order make a lattice with top and bottom.

By omitting the equations $(\hat{k}^{}\bot)$ and
$(\check{k}^{}\top)$ in the definition of $\ml_{\top,\bot}$ we
would obtain the {\it bicartesian} category freely generated by
the set of propositional letters (cf.\ \cite{LS86}, Section I.8).
Dicartesian categories were considered under the name {\it
coherent bicartesian} categories in the printed version of
\cite{DP01a}.

We previously believed wrongly that we have proved coherence for
dicartesian, alias coherent bicartesian, categories. Lemma 5.1 of
the printed version of \cite{DP01a} is however not correct. We
prove here only a restricted coherence result for dicartesian
categories. A study of equality of arrows in bicartesian
categories may be found in \cite{CS01}.

Suppose that in the definition of $\ml_{\top,\bot}$ we omit one of
$\top$ and $\bot$ from the language, and we omit all the arrow
terms and equations involving the omitted nullary connective. When
we omit $\top$, we obtain the category $\ml_{\bot}$, and when we
omit $\bot$, we obtain the category $\ml_{\top}$. It is clear that
$\ml_{\bot}$ is isomorphic to $\ml_{\top}^{op}$. In the printed
version of \cite{DP01a} the categories for which $\ml_{\bot}$ is
the one freely generated by the set of propositional letters were
called {\it coherent sesquicartesian} categories. We call them now
just {\it sesquicartesian} categories.

The category \set, whose objects are sets and whose arrows are
functions, with cartesian product $\times$ as $\kon$, disjoint
union $+$ as $\dis$, a singleton set $\{\ast\}$ as $\top$ and the
empty set $\emptyset$ as $\bot$, is a bicartesian category, but
not a dicartesian category. It is, however, a sesquicartesian
category in the $\ml_{\bot}$ sense, but not in the $\ml_{\top}$
sense. This is because in $\set\,$ we have that
$\emptyset\times\emptyset$ is equal to $\emptyset$, but
$\{\ast\}+\{\ast\}$ is not isomorphic to $\{\ast\}$.

To define the functor $G$ from $\ml_{\top,\bot}$ to \emph{Rel} we
assume that the objects of \emph{Rel} are the formulae of
$\eL_{\top,\bot}$. Everything else in the definition of \emph{Rel}
remains unchanged; in particular, the arrows are sets of ordered
pairs of occurrences of propositional \emph{letters} (no
propositional constant is involved). In the definition of the
functor $G$ we stipulate that for $G\hat{\kappa}_A$ and
$G\check{\kappa}_A$ we have the empty set of ordered pairs. This
serves also for the definition of the functors $G$ from
$\ml_{\bot}$ and $\ml_{\top}$ to \emph{Rel}.

We can establish unrestricted coherence for sesquicartesian
categories, with a proof taken over from the revised version of
\cite{DP01a}, which we will present below. (This proof differs
from the proof in the printed version of \cite{DP01a}, which
relied also on Lemma 5.1, and is not correct.) It is obtained by
enlarging the proof of Lattice Coherence.

The Gentzen formulation of $\ml_{\top,\bot}$ is obtained like that
of \ml\ save that we have in addition the primitive Gentzen terms
$\hat{\kappa}_A\!:A\vdash\top$ and $\check{\kappa}_A\!:\bot\vdash
A$. For Gentzen terms we have as additional equations, besides
\mbox{($\hat{\kappa}$)} and \mbox{($\check{\kappa}$)}, the
following equations:
\begin{tabbing}
\hspace{9em}\=$(\hat{K}^{}\bot)$\hspace{3em}\=$
\hat{K}^{1}_{\bot}\mj_{\bot}$\=$=\hat{K}^{2}_{\bot}\mj_{\bot}$,
\\[1ex]
\>$(\check{K}^{}\top)$\>$\check{K}^{1}_{\top}\mj_{\top}$\>$
=\check{K}^{2}_{\top}\mj_{\top}$,
\end{tabbing}
which amount to $(\hat{k}^{}\bot)$ and $(\check{k}^{}\top)$.

We can prove Composition Elimination for $\ml_{\top,\bot}$ by
enlarging the proof for \ml. We have as new cases first those
where $f$ is $\check{\kappa}_A$ or $g$ is $\hat{\kappa}_A$, which
are taken care of by the equations \mbox{($\check{\kappa}$)} and
\mbox{($\hat{\kappa}$)}. The following case remains. If $f$ is
$\hat{\kappa}_A$, then $g$ is of a form covered by cases already
dealt with. Note that we do not need the equations
$(\hat{K}^{}\bot)$ and $(\check{K}^{}\top)$ for this proof (so
that we have also Composition Elimination for the free bicartesian
category).

Let the category $\mlk_{\dis,\top,\bot}$ be defined like the
category $\mlk_{\dis}$ save that it involves also $\hat{\kappa}$
and the equations \mbox{($\hat{\kappa}$)} and $(\hat{k}^{}\bot)$,
and let the category $\mld_{\kon,\top,\bot}$ be defined like the
category $\mld_{\kon}$ save that it involves also $\check{\kappa}$
and the equations \mbox{($\check{\kappa}$)} and
$(\check{k}^{}\top)$. Composition Elimination is provable for
 $\mlk_{\dis,\top,\bot}$ and
$\mld_{\kon,\top,\bot}$ by abbreviating the proof of Composition
Elimination for $\ml_{\top,\bot}$, in the same way as we
abbreviated the proof of Composition Elimination for \ml\ in order
to obtain Composition Elimination for $\mlk_{\dis}$.

An arrow term of $\ml_{\top,\bot}$ is in {\it standard form} when
it is of the form $g\cirk f$ for $f$ an arrow term
$\mlk_{\dis,\top,\bot}$ and $g$ an arrow term of
$\mld_{\kon,\top,\bot}$. We can then prove the following.

\prop{Standard-Form Lemma}{Every arrow term of $\ml_{\top,\bot}$
is equal in $\ml_{\top,\bot}$ to an arrow term in standard form.}

\dkz By categorial and bifunctorial equations, we may assume that
we deal with a factorized arrow term $f$ none of whose factors is
a complex identity (i.e., $f$ is a big composition of
composition-free arrow terms none of which is equal to an identity
arrow; see \cite{DP04}, Sections 2.6-7, for precise definitions of
these notions) and every factor of $f$ is either an arrow term of
$\mlk_{\dis,\top,\bot}$, and then we call it a {\it
$\kon$-factor}, or an arrow term of $\mld_{\kon,\top,\bot}$, when
we call it a {\it $\dis$-factor}.

Suppose $f\!:B\vdash C$ is a $\kon$-factor and $g\!:A\vdash B$ is
a $\dis$-factor. We show by induction on the length of $f\cirk g$
that in $\ml_{\top,\bot}$
\[
(\ast)\quad f\cirk g=g'\cirk f'\quad{\mbox{\rm or}}\quad f\cirk
g=f'\quad {\mbox{\rm or}}\quad f\cirk g=g'
\]
for $f'$ a $\kon$-factor and $g'$ a $\dis$-factor.

We will consider various cases for $f$. In all such cases, if $g$
is $\check{w}_B$, then we use \mbox{($\check{w}$ {\it nat})}. If
$f$ is $\hat{w}_B$, then we use \mbox{($\hat{w}$ {\it nat})}. If
$f$ is $\hat{k}^{i}_{D,E}$ and $g$ is $g_1\kon g_2$, then we use
\mbox{($\hat{k}^{i}$ {\it nat})}. If $f$ is $f_1\kon f_2$ and $g$
is $g_1\kon g_2$, then we use bifunctorial and categorial
equations and the induction hypothesis.

If $f$ is $f_1\dis f_2$, then we have the following cases. If $g$
is $\check{k}^{i}_{B_1,B_2}$, then we use \mbox{($\check{k}^{i}$
{\it nat})}. If $g$ is $g_1\dis g_2$, then we use bifunctorial and
categorial equations and the induction hypothesis.

Finally, cases where $f$ is $\hat{\kappa}_B$ or $g$ is
$\check{\kappa}_B$ are taken care of by the equations
\mbox{($\hat{\kappa}$)} and \mbox{($\check{\kappa}$)}. This proves
$(\ast)$, and it is clear that $(\ast)$ is sufficient to prove the
lemma. \qed

\vspace{2ex}

We can also prove Composition Elimination and an analogue of the
Standard-Form Lemma for $\ml_{\bot}$. Next we have the following
lemmata for $\ml_{\top,\bot}$ and $\ml_{\bot}$.

\prop{Lemma 1}{If for $f,g\!:A\vdash B$ either $A$ or $B$ is
isomorphic to $\top$ or $\bot$, then $f=g$.}

\dkz If $A$ is isomorphic to $\bot$ or $B$ is isomorphic to
$\top$, then the matter is trivial. Suppose $i\!:B\vdash \bot$ is
an isomorphism. Then from
\[
\hat{k}^{1}_{\bot,\bot}\cirk\langle i\cirk f, i\cirk g\rangle =
\;\hat{k}^{2}_{\bot,\bot}\cirk\langle i\cirk f, i\cirk g\rangle
\]
we obtain $i\cirk f=i\cirk g$, which yields $f=g$. We proceed
analogously if $A$ is isomorphic to $\top$. \qed

\prop{Lemma 2}{If for $f,g\!:A\vdash B$ we have $Gf=Gg=\emptyset$,
then $f=g$.}

\dkz This proof depends on the Standard-Form Lemma above. We write
down $f$ in the standard form $f_2\cirk f_1$ for $f_1\!:A\vdash C$
and $g$ in the standard form $g_2\cirk g_1$ for $g_1\!:A\vdash D$.
Since $\check{k}^{i}$ and $\check{\kappa}$ do not occur in $f_1$,
for every occurrence $z$ of a propositional letter in $C$ we have
an occurrence $x$ of that propositional letter in $A$ such that
$(x,z)\in Gf_1$, and since $\hat{k}^{i}$ and $\hat{\kappa}$ do not
occur in $f_2$, for every occurrence $z$ of a propositional letter
in $C$ we have an occurrence $y$ of that propositional letter in
$B$ such that $(z,y)\in Gf_2$. So if $C$ were not letterless, then
$Gf$ would not be empty. We conclude analogously that $D$, as well
as $C$, is a letterless formula.

If both $C$ and $D$ are isomorphic to $\top$ or $\bot$, then we
have an isomorphism $i\!:C\vdash D$, and $f=f_2\cirk i^{-1}\cirk
i\cirk f_1$. By Lemma~1, we have $i\cirk f_1=g_1$ and $f_2\cirk
i^{-1}=g_2$, from which $f=g$ follows. If $i\!:C\vdash\bot$ and
$j\!:\top\vdash D$ are isomorphisms, then by Lemma ~1 we have
\[
f_2\cirk f_1=g_2\cirk j\cirk\hat{\kappa}_{\bot}\cirk i\cirk f_1
=g_2\cirk g_1,
\]
\noindent and so $f=g$. (Note that
$\hat{\kappa}_{\bot}=\check{\kappa}_{\top}$.) \qed

\vspace{2ex}

We can then prove the following.

\prop{Sesquicartesian Coherence}{The functor $G$ from $\ml_{\bot}$
to \emph{Rel} is faithful.}

\dkz We have Lemma~2 for the case when $Gf=Gg=\emptyset$. When
$Gf=Gg\neq\emptyset$, we proceed as in the proof of Lattice
Coherence, appealing if need there is to Lemma~2, until we reach
the case when $A$ is $A_1\kon A_2$ or a letter, and $B$ is
$B_1\dis B_2$ or a letter, but $A$ and $B$ are not both letters.
In that case, by Composition Elimination, the arrow term $f$ is
equal in $\ml_{\bot}$ either to an arrow term of the form $f'\cirk
\hat{k}^{i}_{A_1,A_2}$, or to an arrow term of the form
$\check{k}^{i}_{B_1,B_2}\cirk f'$. Suppose $f=f'\cirk
\hat{k}^{1}_{A_1,A_2}$. Then for every $(x,y)\in Gf$ we have that
$x$ is in $A_1$. (We reason analogously when $f=f'\cirk
\hat{k}^{2}_{A_1,A_2}$.)

By Composition Elimination too, $g$ is equal in $\ml_{\bot}$
either to an arrow term of the form ${g'\cirk
\hat{k}^{i}_{A_1,A_2}}$, or to an arrow term of the form
${\check{k}^{i}_{B_1,B_2}\cirk g'}$. In the first case we must
have ${g=g'\cirk\hat{k}^{1}_{A_1,A_2}}$, because ${Gg=G(f'\cirk
\hat{k}^{1}_{A_1,A_2})\neq\emptyset}$, and then we apply the
induction hypothesis to derive ${f'=g'}$ from ${Gf'=Gg'}$. Hence
${f=g}$ in $\ml_{\bot}$.

Suppose ${g=\;\check{k}^{1}_{B_1,B_2}\cirk g'}$. (We reason
analogously when ${g=\;\check{k}^{2}_{B_1,B_2}\cirk g'}$.) Let
${f''\!:A_1\vdash B_1\dis B_2''}$ be the substitution instance of
${f'\!:A_1\vdash B_1\dis B_2}$ obtained by replacing every
occurrence of propositional letter in $B_2$ by $\bot$. There is an
isomorphism ${i\!:B_2''\vdash\bot}$, and $f''$ exists because in
$Gf$, which is equal to ${G(\check{k}^{1}_{B_1,B_2}\cirk g')}$,
there is no pair ${(x,y)}$ with $y$ in $B_2$. So we have an arrow
$f'''\!:A_1\vdash B_1$, which we define as
${[\mj_{B_1},\check{\kappa}_{B_1}]\cirk(\mj_{B_1}\dis i)\cirk
f''}$. It is easy to verify that ${G(\check{k}^{1}_{B_1,B_2}\cirk
f''')=Gf'}$, and that ${G(f'''\cirk \hat{k}^{1}_{A_1,A_2})=Gg'}$.
By the induction hypothesis, we obtain
${\check{k}^{1}_{B_1,B_2}\cirk f'''=f'}$ and ${f'''\cirk
\hat{k}^{1}_{A_1,A_2}=g'}$, from which we derive ${f=g}$. We
reason analogously when ${f=\;\check{k}^{i}_{B_1,B_2}\cirk
f'}$.\qed

\vspace{2ex}

From Sesquicartesian Coherence we infer coherence for
$\ml_{\top}$, which is isomorphic to $\ml_{\bot}^{op}$.

\section{\large\bf Restricted coherence for dicartesian categories}
For dicartesian categories we can prove easily a simple restricted
coherence result, which was sufficient for the needs of
\cite{DP04}. A more general, but still restricted, coherence
result with respect to \emph{Rel}, falling short of full
coherence, may be found in the revised version of \cite{DP01a}
(Section 7). We present first the simple restricted coherence
result, and will deal with the more general restricted coherence
result later on.

We define inductively formulae of $\eL_{\top,\bot}$ in {\it
disjunctive normal form} ({\it dnf}$\,$): every $\dis$-free
formula is in {\it dnf}, and if $A$ and $B$ are both in {\it dnf},
then $A\dis B$ is in {\it dnf}. We define dually formulae of
$\eL_{\top,\bot}$ in {\it conjunctive normal form} ({\it
cnf}$\,$): every $\kon$-free formula is in {\it cnf}, and if $A$
and $B$ are both in {\it cnf}, then $A\kon B$ is in {\it cnf}.

\prop{Restricted Dicartesian Coherence}{Let $f,g\!:A\vdash B$ be
arrow terms of $\ml_{\top,\bot}$ such that $A$ is in {\rm dnf} and
$B$ in {\rm cnf}. If $Gf=Gg$, then $f=g$ in $\ml_{\top,\bot}$.}

\dkz If ${Gf=Gg=\emptyset}$, then we apply Lemma 2. If
${Gf=Gg\neq\emptyset}$, then we proceed as in the proof of Lattice
Coherence, by induction on the sum of the lengths of $A$ and $B$,
appealing if need there is to Lemma~2, until we reach the case
when $A$ is $A_1\kon A_2$ or a letter, and $B$ is $B_1\dis B_2$ or
a letter, but $A$ and $B$ are not both letters. In that case there
is no occurrence of $\dis$ in $A$ and no occurrence of $\kon$ in
$B$. We then rely on the composition-free form of $f$ and $g$ in
$\ml_{\top,\bot}$ and on the equation $(\hat{K}\check{K})$.\qed

\vspace{2ex}

To improve upon this result we need the following lemma for
$\ml_{\top,\bot}$, and the definitions that follow. This lemma is
analogous up to a point to the Invertibility Lemma for $\dis$.

\prop{Lemma 3}{Let ${f: A\vdash B_1\dis B_2}$ be a Gentzen term
such that ${Gf\neq\emptyset}$ and $\dis$ does not occur in $A$. If
for every ${(x,y)}\in Gf$ we have that $y$ is in $B_1$, then there
is a Gentzen term ${g: A\vdash B_1}$ such that
${Gf=G\check{K}^1_{B_2}g}$.}

\dkz We proceed by induction on the length of $A$. Suppose $f$ is
a composition-free Gentzen term. If $A$ is a propositional letter,
then by the assumption on ${Gf}$ we have that $f$ is of the form
${\check{K}^1_{B_2}f'}$, and we can take that $g$ is $f'$.

If $A$ is not a propositional letter and $f$ is not of the form
${\check{K}^1_{B_2}f'}$ (by the assumption on ${Gf}$, the Gentzen
term $f$ cannot be of the form ${\check{K}^2_{B_1}f'}$), then,
since $\dis$ does not occur in $A$, we have that $f$ is of the
form ${\hat{K}^i_{A''}f'}$ for ${f':A'\vdash B_1\dis B_2}$. Note
that ${Gf'\neq\emptyset}$ and $\dis$ does not occur in $A'$. Since
for every ${(x,y)}$ in $Gf'$ we have that $y$ is in $B_1$, we may
apply the induction hypothesis to $f'$ and obtain $g'\!:A'\vdash
B_1$ such that ${Gf'=G\check{K}^1_{B_2}g'}$. By relying on the
equation $(\hat{K}\check{K})$, we can take that $g$ is
${\hat{K}^i_{A''}g'}$. \qed

\vspace{2ex}

A formula $C$ of $\eL_{\top,\bot}$ is called a {\em contradiction}
when there is in $\ml_{\top,\bot}$ an arrow of the type ${C\vdash
\bot}$. For every formula that is not a contradiction there is a
substitution instance isomorphic to $\top$. Suppose $C$ is not a
contradiction, and let $C^{\top}$ be obtained from $C$ by
substituting $\top$ for every propositional letter. If $C^{\top}$
were not isomorphic to $\top$, then since every letterless formula
of $\eL_{\top,\bot}$ is isomorphic in $\ml_{\top,\bot}$ either to
$\top$ or to $\bot$, we would have an isomorphism
${i:C^{\top}\vdash \bot}$. Since there is obviously an arrow
${u:C\vdash C^{\top}}$ formed by using $\hat{\kappa}_p$, we would
have ${i\cirk u:C\vdash\bot}$, and $C$ would be a contradiction.

A formula $C$ of $\eL_{\top,\bot}$ is called a {\em tautology}
when there is in $\ml_{\top,\bot}$ an arrow of the type
$\top\vdash C$. For every formula that is not a tautology there is
a substitution instance isomorphic to $\bot$. (This is shown
analogously to what we had in the preceding paragraph.)

A formula of $\eL_{\top,\bot}$ is called $\bot$-{\em normal} when
for every subformula $D\kon C$ or $C\kon D$ of it with $C$ a
contradiction, there is no occurrence of $\dis$ in $D$. A formula
of $\eL_{\top,\bot}$ is called $\top$-{\em normal} when for every
subformula $D\dis C$ or $C\dis D$ of it with $C$ a tautology,
there is no occurrence of $\kon$ in $D$.

We can now formulate our second partial coherence result for
dicartesian categories.

\prop{Restricted Dicartesian Coherence II}{If ${f,g:A\vdash B}$
are terms of $\ml_{\top,\bot}$ such that ${Gf=Gg}$ and either $A$
is $\bot$-normal or $B$ is $\top$-normal, then ${f=g}$ in
$\ml_{\top,\bot}$.}

\dkz Suppose $A$ is $\bot$-normal. Lemma~2 covers the case when
$Gf=Gg=\emptyset$. So we assume $Gf=Gg\neq\emptyset$, and proceed
as in the proof of Sesquicartesian Coherence by induction on the
sum of the lengths of $A$ and $B$. The basis of this induction and
the cases when $A$ is of the form ${A_1\dis A_2}$ or $B$ is of the
form ${B_1\kon B_2}$ are settled as in the proof of
Sesquicartesian Coherence.

Suppose $A$ is ${A_1\kon A_2}$ or a propositional letter and $B$
is ${B_1\dis B_2}$ or a propositional letter, but $A$ and $B$ are
not both propositional letters. (The cases when $A$ or $B$ is a
constant object are excluded by the assumption that
$Gf=Gg\neq\emptyset$.) We proceed then as in the proof of
Sesquicartesian Coherence until we reach the case when ${f=f'\cirk
\hat{k}^1_{A_1,A_2}}$ and ${g=\check{k}^1_{B_1,B_2}\cirk g'}$.

Suppose $A_2$ is not a contradiction. Then there is an instance
${A_2^{\top}}$ of $A_2$ and an isomorphism ${i:\top\vdash
A_2^{\top}}$. (To obtain ${A_2^{\top}}$ we substitute $\top$ for
every letter in $A_2$.) Let ${g'':A_1\kon A_2^{\top}\vdash B_1}$
be the substitution instance of ${g':A_1\kon A_2\vdash B_1}$
obtained by replacing every occurrence of propositional letter in
$A_2$ by $\top$. Such a term exists because in ${Gg}$, which is
equal to ${G(f'\cirk \hat{k}^1_{A_1,A_2})}$, there is no pair
${(x,y)}$ with $x$ in $A_2$.

So we have an arrow $g'''=g''\cirk(\mj_{A_1}\kon
i)\cirk\langle\mj_{A_1},\hat{\kappa}_{A_1}\rangle:A_1\vdash B_1$.
It is easy to verify that $G(\check{k}^1_{B_1,B_2}\cirk g''')=Gf'$
and that $G(g'''\cirk \hat{k}^1_{A_1,A_2})=Gg'$. By the induction
hypothesis we obtain $\check{k}^1_{B_1,B_2}\cirk g'''=f'$ and
$g'''\cirk \hat{k}^1_{A_1,A_2}=g'$, from which we derive $f=g$.

Suppose $A_2$ is a contradiction. Then by the assumption that $A$
is $\bot$-normal we have that $\dis$ does not occur in $A_1$. We
may apply Lemma~3 to ${f':A_1\vdash B_1\dis B_2}$ to obtain
$f''':A_1\vdash B_1$ such that $Gf'=G(\check{k}^1_{B_1,B_2}\cirk
f''')$. It is easy to verify that then $Gg'=G(f'''\cirk
\hat{k}^1_{A_1,A_2})$, and we may proceed as in the proof of
Sesquicartesian Coherence.

We proceed analogously when $B$ is $\top$-normal, relying on a
lemma dual to Lemma~3.\qed

\vspace{2ex}

Consider the following definitions:
\begin{tabbing}
\hspace{9em}\=$A^0_{\bot}$\=$=A\kon
\bot$,\hspace{2em}\=$A^{n+1}_{\bot}$\=$=(A^n_{\bot}\dis\top)\kon
\bot$,\\*[1ex]

\>$f^0_{\bot}$\>$=f\kon \mj_{\bot}$,\>
$f^{n+1}_{\bot}$\>$=(f^n_{\bot}\dis\mj_\top)\kon
\mj_{\bot}$,\\[2ex]\>$A^0_\top$\>$=A\dis \top$,\>$A^{n+1}_\top$\>
$=(A^n_\top\kon\bot)\dis \top$,\\*[1ex]\>$f^0_\top$\>$=f\dis
\mj_\top$,\>$f^{n+1}_\top$\>$=(f^n_\top\kon\mj_{\bot})\dis
\mj_\top$.
\end{tabbing}
Then for $f^n$ being

\[
(\check{k}^1_{A,\top}\kon\mj_{\bot})^n_\top\cirk \hat{k}^1_{(A\kon
{\bot})^n_\top, {\bot}}: A^{n+1}_{\bot}\vdash A^{n+1}_\top
\]

\noindent and $g^n$ being

\[
\check{k}^1_{(A\dis \top)^n_\bot, \top}\cirk
(\hat{k}^1_{A,{\bot}}\dis\mj_\top)^n_{\bot}: A^{n+1}_{\bot}\vdash
A^{n+1}_\top
\]

\noindent we have $Gf^n=Gg^n$, but we suppose that $f^n=g^n$ does
not hold in $\ml_{\top,\bot}$. The equation $f^0=g^0$ is

\begin{tabbing}

\quad$((\check{k}^1_{A,\top}\kon\mj_{\bot})\dis\mj_\top)\cirk
\hat{k}^1_{(A\kon {\bot})\dis\top, {\bot}}=\check{k}^1_{(A\dis
\top)\kon\bot, \top}\cirk
(\hat{k}^1_{A,{\bot}}\dis\mj_\top)\kon\mj_{\bot}:$
\\[.5ex]
\`$((A\kon\bot)\dis\top)\kon\bot\vdash
((A\dis\top)\kon\bot)\dis\top.$
\end{tabbing}

\noindent Note that $A^{n+1}_{\bot}$ is not $\bot$-normal, and
$A^{n+1}_\top$ is not $\top$-normal.

We don't know whether it is sufficient to add to $\ml_{\top,\bot}$
the equations $f^n=g^n$ for every $n\geq 0$ in order to obtain
full coherence for the resulting category.

As a corollary of Restricted Dicartesian Coherence II, we obtain
that if ${f,g:A\vdash B}$ are terms of $\ml_{\top,\bot}$ such that
${Gf=Gg}$, while $A$ and $B$ are isomorphic either to formulae of
\eL\ (i.e.\ to formulae in which $\top$ and $\bot$ do not occur)
or to letterless formulae, then ${f=g}$ in $\ml_{\top,\bot}$. This
corollary is analogous to the restricted coherence result for
symmetric monoidal closed categories of Kelly and Mac Lane in
\cite{KML71} (see \cite{DP07}, Section 3.1).

\section{\large\bf Maximality}
A syntactically built category such as \ml\ and $\ml_{\top,\bot}$
is called {\it maximal} when adding any new axiomatic equation
between arrow terms of this category yields a category that is a
preorder. The new axiomatic equation is supposed to be closed
under substitution for propositional letters, as the equations of
\ml\ and $\ml_{\top,\bot}$ were. (This notion of maximality for
syntactical categories is defined more precisely in \cite{DP04},
Section 9.3.) Maximality is an interesting property when the
initial category, like \ml\ and $\ml_{\top,\bot}$ here, is not
itself a preorder. We will deal in subsequent sections with
maximality for \ml\ and $\ml_{\top,\bot}$.

The maximality property above is analogous to the property of
usual formulations of the classical propositional calculus called
{\it Post completeness}. That this calculus is Post complete means
that if we add to it any new axiom-schema in the language of the
calculus, then we can prove every formula. An analogue of B\"
ohm's Theorem in the typed lambda calculus implies, similarly,
that the typed lambda calculus cannot be extended without falling
into triviality, i.e.\ without every equation (between terms of
the same type) becoming derivable (see \cite{S95}, \cite{DP00} and
references therein; see \cite{B81}, Section 10.4, for B\" ohm's
Theorem in the untyped lambda calculus).

Let us now consider several examples of common algebraic
structures with analogous maximality properties. First, we have
that semilattices are maximal in the following sense.

Let $a$ and $b$ be terms made exclusively of variables and of a
binary operation $\cdot$, which we interpret as meet or join. That
the equation $a=b$ holds in a semilattice $S$ means that {\it
every} instance of $a=b$ obtained by substituting names of
elements of $S$ for variables holds in $S$. Suppose $a=b$ does not
hold in a free semilattice $S_F$ (so it is not the case that $a=b$
holds in every semilattice). Hence there must be an instance of
$a=b$ obtained by substituting names of elements of $S_F$ for
variables such that this instance does not hold in $S_F$. It is
easy to conclude that in $a=b$ there must be at least two
variables, and that $S_F$ must have at least two free generators.
Then every semilattice in which $a=b$ holds is trivial---namely,
it has a single element.

Here is a short proof of that. If $a=b$ does not hold in $S_F$,
then there must be a variable $x$ in one of $a$ and $b$ that is
not in the other. Then from $a=b$, by substituting $y$ for every
variable in $a$ and $b$ different from $x$, and by applying the
semilattice equations, we infer either $x=y$ or $x\cdot y=y$. If
we have $x=y$, we are done, and, if we have $x\cdot y=y$, then we
have also $y\cdot x=x$, and hence $x=y$.

Semilattices with unit, distributive lattices, distributive
lattices with top and bottom, and Boolean algebras are maximal in
the same sense. The equations $a=b$ in question are equations
between terms made exclusively of variables and the operations of
the kind of algebra we envisage: semilattices with unit,
distributive lattices, etc. That such an equation holds in a
particular structure means, as above, that every substitution
instance of it holds. However, the number of variables in $a=b$
and the number of generators of the free structure mentioned need
not always be at least two.

If we deal with semilattices with unit \mj, then $a=b$ must have
at least one variable, and the free semilattice with unit must
have at least one free generator. We substitute \mj\ for every
variable in $a$ and $b$ different from $x$ in order to obtain
$x=\mj$, and hence triviality. So semilattices with unit are
maximal in the same sense.

The same sort of maximality can be proven for distributive
lattices, whose operations are $\kon$ and $\dis$, which we call
conjunction and disjunction, respectively. Then every term made of
$\kon$, $\dis$ and variables is equal to a term in disjunctive
normal form (i.e.\ a multiple disjunction of multiple conjunctions
of variables; see the preceding section for a precise definition),
and to a term in conjunctive normal form (i.e.\ a multiple
conjunction of multiple disjunctions of variables; see the
preceding section). These normal forms are not unique. If $a=b$,
in which we must have at least two variables, does not hold in a
free distributive lattice $D_F$ with at least two free generators,
then either $a\leq b$ or $b\leq a$ does not hold in $D_F$. Suppose
$a\leq b$ does not hold in $D_F$. Let $a'$ be a disjunctive normal
form of $a$, and let $b'$ be a conjunctive normal form of $b$. So
$a'\leq b'$ does not hold in $D_F$. From that we infer that for a
disjunct $a''$ of $a'$ and for a conjunct $b''$ of $b'$ we do not
have $a''\leq b''$ in $D_F$. This means that there is no variable
in common in $a''$ and $b''$; otherwise, the conjunction of
variables $a''$ would be lesser than or equal in $D_F$ to the
disjunction of variables $b''$. If in a distributive lattice $a=b$
holds, then $a''\leq b''$ holds too, and hence, by substitution,
we obtain $x\leq y$. So $x=y$.

For distributive lattices with top $\top$ and bottom $\bot$, we
proceed analogously via disjunctive and conjunctive normal form.
Here $a=b$ may be even without variables, and the free structure
may have even an empty set of free generators. The additional
cases to consider are when in $a''\leq b''$ we have that $a''$ is
$\top$ and $b''$ is $\bot$. In any case, we obtain $\top\leq\bot$,
and hence our structure is trivial.

The same sort of maximality can be proven for Boolean algebras,
i.e.\ complemented distributive lattices. Boolean algebras must
have top and bottom. In a disjunctive normal form now the
disjuncts are conjunctions of variables $x$ or terms $\bar x$,
where $\;\bar{}\;$ is complementation, or the disjunctive normal
form is just $\top$ or $\bot$; analogously for conjunctive normal
forms. Then we proceed as for distributive lattices with an
equation $a=b$ that may be even without variables, until we reach
that $a''\leq b''$, which does not hold in a free Boolean algebra
$B_F$, whose set of free generators may be even empty, holds in
our Boolean algebra. If $x$ is a conjunct of $a''$, then in $b''$
we cannot have a disjunct $x$; but we may have a disjunct $\bar
x$. The same holds for the conjuncts $\bar x$ of $a''$. It is
excluded that both $x$ and $\bar x$ are conjuncts of $a''$, or
disjuncts of $b''$; otherwise, $a''\leq b''$ would hold in $B_F$.
Then for every conjunct $x$ in $a''$ and every disjunct $\bar y$
in $b''$ we substitute $\top$ for $x$ and $y$, and for every other
variable we substitute $\bot$. In any case, we obtain
$\top\leq\bot$, and hence our Boolean algebra is trivial. This is
essentially the proof of Post completeness for the classical
propositional calculus, due to Bernays and Hilbert (see
\cite{Z99}, Section 2.4, and \cite{HA28}, Section I.13), from
which we can infer the ordinary completeness of this calculus with
respect to valuations in the two-element Boolean algebra---namely,
with respect to truth tables---and also completeness with respect
to any nontrivial model.

As examples of common algebraic structures that are not maximal in
the sense above, we have semigroups, commutative semigroups,
lattices, and many others. What is maximal for semilattices and is
not maximal for lattices is the equational theory of the
structures in question. The equational theory of semilattices
cannot be extended without falling into triviality, while the
equational theory of lattices can be extended with the
distributive law, for example.

The notions of maximality envisaged in this section were extreme
(or should we say ``maximal''), in the sense that we envisaged
collapsing only into preorder. For semilattices, distributive
lattices, etc., this is also preorder for a one-object category.
We may, however, envisage relativizing our notion of maximality by
replacing preorder with a weaker property, such that structures
possessing it are trivial, but not so trivial (cf. \cite{D99},
Section 4.11). We will encounter maximality in such a relative
sense in the last section.

As an example of relative maximality in a common algebraic
structure we can take symmetric groups. Consider the standard
axioms for the symmetric group $\eS_n$, where $n\geq 2$, with the
generators $s_i$, for $i\in\{1,\ldots,n\mn 1\}$, corresponding to
transpositions of immediate neighbours (see \cite{CM57}, Section
6.2). If to $\eS_n$ for $n\geq 5$ we add an equation $a=\mj$ where
$a$ is built exclusively of the generators $s_i$ of $\eS_n$ with
composition, and $a=\mj$ does not hold in $\eS_n$, then we can
derive $s_i=s_j$. This does not mean that the resulting structure
will be a one-element structure, i.e.\ the trivial one-element
group. It will be such if $a$ is an odd permutation, and if $a$ is
an even permutation, then we will obtain a two-element structure,
which is $\eS_2$. This can be inferred from facts about the normal
subgroups of $\eS_n$. Simple groups are maximal in the nonrelative
sense, envisaged above for semilattices.

\section{\large\bf Maximality of lattice categories}
We will show in this section that \ml\ is maximal in the sense
specified at the beginning of the preceding section; namely, in
the interesting way. (We take over this result from \cite{DP02},
Section~5, and \cite{DP04}, Section 9.5.)

Suppose $A$ and $B$ are formulae of $\eL$ in which only $p$ occurs
as a letter. If for some arrow terms $f_1,f_2\!: A\vdash B$ of
\ml\ we have $Gf_1\neq Gf_2$, then for some $x$ in $A$ and some
$y$ in $B$ we have $(x,y)\in Gf_1$ and $(x,y)\not\in Gf_2$, or
vice versa. Suppose $(x,y)\in Gf_1$ and $(x,y)\not\in Gf_2$.

For every subformula $C$ of $A$ and every formula $D$ let $A^C_D$
be the formula obtained from $A$ by replacing the particular
occurrence of the formula $C$ in $A$ by $D$. It can be shown that
for every subformula $A_1\dis A_2$ of $A$ we have an arrow term
$h\!:A^{A_1\dis A_2}_{A_j}\vdash A$ of \ml, built by using
$\check{k}^{j}_{A_1,A_2}$, such that there is an $x'$ in
$A^{A_1\dis A_2}_{A_j}$ for which $(x',x)\in Gh$. Hence, for such
an $h$, we have $(x',y)\in G(f_1\cirk h)$ and $(x',y)\not\in
G(f_2\cirk h)$.

We compose $f_i$ repeatedly with such  arrow terms until we obtain
the arrow terms $f'_i\!:p\kon\ldots\kon p\vdash B$ of \ml\ such
that parentheses are somehow associated in $p\kon\ldots\kon p$ and
for some $z$ in $(p\kon\ldots\kon p)$ we have $(z,y)\in Gf'_1$ and
$(z,y)\not\in Gf'_2$. The formula $p\kon\ldots\kon p$ may also be
only $p$. We may further compose $f'_i$ with other arrow terms of
\ml\ in order to obtain the arrow terms $f''_i$ of type $p\kon
A'\vdash B$ or $p\vdash B$ such that $A'$ is of the form
$p\kon\ldots\kon p$ with parentheses somehow associated. Let us
use $0$ to denote the first occurrence of a propositional letter
in a formula, counting from the left. So we have $(0,y)\in Gf''_1$
but $(0,y)\not\in Gf''_2$.

By working dually on $B$ we obtain the arrow terms $f'''_i$ of
\ml\ of type $p\kon A'\vdash p\dis B'$, for $A'$ of the form
$p\kon\ldots\kon p$ and $B'$ of the form $p\dis\ldots\dis p$, or
of type $p\kon A'\vdash p$, or of type $p\vdash p\dis B'$, such
that $(0,0)\in Gf'''_1$ and $(0,0)\not\in Gf'''_2$. (We cannot
obtain that $f'''_1$ and $f'''_2$ are of type $p\vdash p$, since,
otherwise, by Composition Elimination for \ml, $f'''_2$ would not
exist.)

There is an arrow term $h^{\kon}\!:p\vdash p\kon\ldots\kon p$ of
\ml\ defined by using $\hat{w}$ such that for every $x\in
G(p\kon\ldots\kon p)$ we have $(0,x)\in Gh^{\kon}$. We define
analogously with the help of $\check{w}$ an arrow term
$h^{\dis}\!:p\dis\ldots\dis p\vdash p$ of \ml\ such that for every
$x$ in $p\dis\ldots\dis p$ we have $(x,0)\in Gh^{\dis}$. The arrow
terms $h^{\kon}$ and $h^{\dis}$ may be $\mj_p\!:p\vdash p$.

If $f'''_i$ is of type $p\kon A'\vdash p\dis B'$, let
$f^{\dag}_i\!:p\kon p\vdash p\dis p$ be defined by
\[
f^{\dag}_i=_{\df}\; (\mj_p\dis h^{\dis})\cirk
f'''_i\cirk(\mj_p\dis h^{\kon}).
\]
By Composition Elimination for \ml, we have that $Gf^{\dag}_i$
must be a singleton. Let us use $1$ to denote the second
occurrence of a propositional letter in a formula, counting from
the left. If $(1,0)$ or $(1,1)$ belongs to $Gf^{\dag}_2$, then for
$f^{\ast}_i\!:p\kon p\vdash p$ defined as $\check{w}_p\cirk
f^{\dag}_i$ we have $(0,0)\in Gf^{\ast}_1$ and $(0,0)\not\in
Gf^{\ast}_2$. If $(0,1)$ or $(1,1)$ belongs to $Gf^{\dag}_2$, then
for $f^{\ast}_i\!:p\vdash p\dis p$ defined as
$f^{\dag}_i\cirk\hat{w}_p$ we have $(0,0)\in Gf^{\ast}_1$ and
$(0,0)\not\in Gf^{\ast}_2$.

If $f'''_i$ is of type $p\kon A'\vdash p$, then for
$f^{\ast}_i\!:p\kon p\vdash p$ defined as $f'''_i\cirk(\mj_p\dis
h^{\kon})$ we have $(0,0)\in Gf^{\ast}_1$ and $(0,0)\not\in
Gf^{\ast}_2$.

If $f'''_i$ is of type $p\vdash p\dis B'$, then for
$f^{\ast}_i\!:p\vdash p\dis p$ defined as $(\mj_p\dis
h^{\dis})\cirk f'''_i$ we have $(0,0)\in Gf^{\ast}_1$ and
$(0,0)\not\in Gf^{\ast}_2$. In all that we have by Composition
Elimination for \ml\ that $Gf^{\ast}_i$ must be a singleton.

In cases where $f^{\ast}_i$ is of type $p\kon p\vdash p$, by
Composition Elimination for \ml, by the conditions on
$Gf^{\ast}_1$ and $Gf^{\ast}_2$, and by the functoriality of $G$,
we obtain in \ml\ the equation $f^{\ast}_i=\:\hat{k}^{i}_{p,p}$.
(This follows from Lattice Coherence too.) So in \ml\ extended
with $f_1=f_2$ we can derive the equation
\[
(\hat{k}^{}\hat{k}^{})\quad
\hat{k}^{1}_{p,p}\;=\;\hat{k}^{2}_{p,p}.
\]

In cases where $f^{\ast}_i$ is of type $p\vdash p\dis p$, we
conclude analogously that we have in \ml\ the equation
$f^{\ast}_i=\:\check{k}^{i}_{p,p}$, and so in \ml\ extended with
$f_1=f_2$ we can derive
\[
(\check{k}^{}\check{k}^{})\quad
\check{k}^{1}_{p,p}\;=\;\check{k}^{2}_{p,p}.
\]
If either of $(\hat{k}^{}\hat{k}^{})$ and
$(\check{k}^{}\check{k}^{})$ holds in a lattice category \aA, then
\aA\ is a preorder.

It remains to remark that if for some arrow terms $g_1$ and $g_2$
of \ml\ of the same type we have that $g_1=g_2$ does not hold for
\ml, then by Lattice Coherence we have $Gg_1\neq Gg_2$. If we take
the substitution instances $g'_1$ of $g_1$ and $g'_2$ of $g_2$
obtained by replacing every letter by a single letter $p$, then we
obtain again $Gg'_1\neq Gg'_2$. If $g_1=g_2$ holds in a lattice
category \aA, then $g'_1=g'_2$ holds too, and \aA\ is a preorder,
as we have shown above. This concludes the proof of maximality for
\ml. (In the original presentation of this proof in \cite{DP02},
Section~5, there are some slight inaccuracies in the definition of
$f^{\ast}_i$.)

\section{\large\bf Relative maximality of dicartesian categories}
The category $\ml_{\top,\bot}$ is not maximal in the sense in
which \ml\ is. This is shown by the following counterexample.

Let $\set_{\ast}$ be the category whose objects are sets with a
distinguished element $\ast$, and whose arrows are
$\ast$-preserving functions $f$ between these sets; namely,
$f(\ast)=\ast$. This category is isomorphic to the category of
sets with partial functions. The following definitions serve to
show that $\set_{\ast}$ is a category in which we can interpret
the objects and arrow terms of $\ml_{\top,\bot}$:
\begin{tabbing}
\quad\quad I$\;=\{\ast\}$, \quad\quad \= $a'=\{(x,\ast)\mid x\in
a-{\mbox{\rm I}\}}$, \quad\quad $b''=\{(\ast,y)\mid y\in
b-{\mbox{\rm I}}\}$,
\\[1.5ex]
\> $a\otimes b\:$ \= = \= $((a-{\mbox{\rm I}})\kon(b-{\mbox{\rm
I}}))\cup{\mbox{\rm I}}$,
\\
\> $a\Bx b\:$ \> = \> $(a\otimes b)\cup a'\cup b''$,
\\
\> $a\Bk b\:$ \> = \> $a'\cup b''\cup {\mbox{\rm I}}$.
\end{tabbing}
Note that $a\Bx b$ is isomorphic in $\set\,$ to the cartesian
product $a\times b$; the element $\ast$ of $a\Bx b$ corresponds to
the element $(\ast,\ast)$ of $a\times b$.

The functions $\hat{k}^{i}_{a_1,a_2}:a_1\Bx a_2\str a_i$, for
$i\in\{1,2\}$, are defined by
\[
\hat{k}^{i}_{a_1,a_2}\!(x_1,x_2)=x_i, \quad\quad
\hat{k}^{i}_{a_1,a_2}\!(\ast)=\ast;
\]
for $f_i\!:c\str a_i$, the function $\langle f_1,f_2\rangle\!:
c\str a_1\Bx a_2$ is defined by
\[
\langle f_1,f_2\rangle(z)= \left\{
\begin{array}{ll}
(f_1(z),f_2(z)) & {\mbox{\rm if }} f_1(z)\neq\ast {\mbox{\rm{ or
}}} f_2(z)\neq\ast
\\[.5ex]
\ast & {\mbox{\rm if }} f_1(z)=f_2(z)=\ast;
\end{array}
\right.
\]
and the function $\hat{\kappa}_a:a\str {\mbox{\rm I}}$ is defined
by $\hat{\kappa}_a(x)=\ast$. Having in mind the isomorphism
between $a\Bx b$ and $a\times b$ mentioned above, the functions
$\hat{k}^{i}_{a_1,a_2}:a_1\Bx a_2\str a_i$ correspond to the
projection functions, while \mbox{$\langle\_\,,\_\,\rangle$}
corresponds to the usual pairing operation on functions.

The functions $\check{k}^{i}_{a_1,a_2}:a_i\str a_1\Bk a_2$ are
defined by
\[
\begin{array}{l}
\check{k}^{1}_{a_1,a_2}(x)=(x,\ast),\quad
\check{k}^{2}_{a_1,a_2}(x)=(\ast,x), \quad{\mbox{\rm{for }}}
x\neq\ast,
\\[.5ex]
\check{k}^{i}_{a_1,a_2}(\ast)=\ast;
\end{array}
\]
for $g_i\!:a_i\str c$, the function $[g_1,g_2]\!:a_1\Bk a_2 \str
c$ is defined by
\[
\begin{array}{l}
[g_1,g_2](x_1,x_2)=g_i(x_i), {\mbox{\rm{ for }}} x_i\neq\ast,
\\[.5ex]
[g_1,g_2](\ast)=\ast;
\end{array}
\]
finally, the function $\check{\kappa}_a:{\mbox{\rm I}}\str a$ is
defined by $\check{\kappa}_a(\ast)=\ast$.

If we take that $\kon$ is $\Bx$ and $\dis$ is $\Bk\!\!$, then it
can be checked in a straightforward manner that $\set_{\ast}$ and
$\set_{\ast}$ without {\mbox{\rm I}} are lattice categories, and
if in $\set_{\ast}$ we take further that both $\top$ and $\bot$
are I, then $\set_{\ast}$ is a dicartesian category.

Consider now the category $\set_{\ast}^{\emptyset}$, which is
obtained by adding to $\set_{\ast}$ the empty set $\emptyset$ as a
new object, and the empty functions $\emptyset_a\!:\emptyset\str
a$ as new arrows. The identity arrow $\mj_{\emptyset}$ is
$\emptyset_{\emptyset}$. For $\set_{\ast}^{\emptyset}$, we enlarge
the definitions above by
\[
\begin{array}{rl}
\emptyset\Bx a\!\!\!&=a\Bx\emptyset=\emptyset,\\
\emptyset\Bk a\!\!\!&=a\Bk\emptyset=a,\\[1ex]
\hat{k}^{i}_{a_1,a_2}\!\!\!&=\emptyset_{a_i}, {\mbox{\rm{ for }}}
a_1=\emptyset {\mbox{\rm{ or }}} a_2=\emptyset,\\[.5ex]
\langle \emptyset_{a_1},\emptyset_{a_2}\rangle\!\!\!&=
\emptyset_{a_1\!\!\bx\,a_2},
\\[.5ex]
\hat{\kappa}_{\emptyset}\!\!\!&=\emptyset_{\mbox{\rm{\scriptsize I}}},\\[1ex]
\check{k}^{i}_{a_1,a_2}\!\!\!&=\emptyset_{a_1\!\!\bk\,a_2},
{\mbox{\rm{ for }}}
a_i=\emptyset,\\[.5ex]
[f_1,\emptyset_c]\!\!\!&=f_1,\quad\quad\quad
[\emptyset_c,f_2]=f_2,
\end{array}
\]
and define now the function $\check{\kappa}_a:\emptyset\str a$ by
$\check{\kappa}_a=\emptyset_a$. Then it can be checked that
$\set_{\ast}^{\emptyset}$ where $\kon$ is $\Bx$ and $\dis$ is
$\Bk$ as before, while $\top$ is I and $\bot$ is $\emptyset$, is a
dicartesian category too.

In $\ml_{\top,\bot}$ the equation
$\hat{k}^{1}_{p,\bot}=\check{\kappa}_p\cirk\hat{k}^{2}_{p,\bot}$
does not hold, because $G\hat{k}^{1}_{p,\bot}\neq\emptyset$ and
$G(\check{\kappa}_p\cirk\hat{k}^{2}_{p,\bot})=\emptyset$, but in
$\set_{\ast}^{\emptyset}$ this equation holds, because both sides
are equal to $\emptyset_{\emptyset}$. Since
$\set_{\ast}^{\emptyset}$ is not a preorder, we can conclude that
$\ml_{\top,\bot}$ is not maximal.

Although this maximality fails, the category $\ml_{\top,\bot}$ may
be shown maximal in a relative sense. This relative maximality
result, which we are going to demonstrate now, says that every
dicartesian category that satisfies an equation $f=g$ between
arrow terms of $\ml_{\top,\bot}$ such that $Gf\neq Gg$ (which
implies that $f=g$ is not in $\ml_{\top,\bot}$) satisfies also
some particular equations. These equations do not give preorder in
general, but a kind of ``contextual'' preorder. Moreover, when
$\ml_{\top,\bot}$ is extended with some of these equations we
obtain a maximal category.

If for some arrow terms $f_1,f_2\!: A\vdash B$ of
$\ml_{\top,\bot}$ we have $Gf_1\neq Gf_2$, then for some $x$ in
$A$ and some $y$ in $B$ we have $(x,y)\in Gf_1$ and $(x,y)\not\in
Gf_2$, or vice versa. Suppose $(x,y)\in Gf_1$ and $(x,y)\not\in
Gf_2$. Suppose $x$ is an occurrence of $p$, so that $y$ must be an
occurrence of $p$ too.

Let $A'$ be the formula obtained from the formula $A$ by replacing
$x$ by $p\kon\bot$, and every other occurrence of letter or $\top$
by $\bot$. Dually, let $B'$ be the formula obtained from $B$ by
replacing $y$ by $p\dis\top$, and every other occurrence of letter
or $\bot$ by $\top$. Let us use $0$, as in the preceding section,
to denote the first occurrence of a propositional letter in a
formula, counting from the left. Then it can be shown that there
is an arrow term $h^A\!:A'\vdash A$ of $\ml_{\top,\bot}$ such that
$Gh^A=\{(0,x)\}$, and an arrow term $h^B\!:B\vdash B'$ of
$\ml_{\top,\bot}$ such that $Gh^B=\{(y,0)\}$. We build $h^A$ with
$\hat{k}^{1}_{p,\bot}:p\kon\bot\vdash p$ and instances of
$\check{\kappa}_C:\bot\vdash C$, with the help of the operations
$\kon$ and $\dis$ on arrow terms. Analogously, $h^B$ is built with
$\check{k}^{1}_{p,\top}:p\vdash p\dis\top$ and instances of
$\hat{\kappa}_C:C\vdash\top$. It can also be shown that there are
arrow terms $j^A\!:p\kon\bot\vdash A'$ and $j^B\!:B'\vdash
p\dis\top$ of $\ml_{\top,\bot}$ such that $Gj^A=Gj^B=\{(0,0)\}$.
These arrow terms stand for isomorphisms of $\ml_{\top,\bot}$.

Then it is clear that for $f'_i$ being
\[
j^B\cirk h^B\cirk f_i\cirk h^A\cirk j^A\!:p\kon\bot\vdash
p\dis\top,
\]
with $i\in\{1,2\}$, we have $Gf'_1=\{(0,0)\}$, while
$Gf'_2=\emptyset$. Hence, by Composition Elimination for
$\ml_{\top,\bot}$ and by the functoriality of $G$, we obtain in
$\ml_{\top,\bot}$ the equations
\[
\begin{array}{l}
f'_1=\;\check{k}^{1}_{p,\top}\cirk\hat{k}^{1}_{p,\bot},\\[.5ex]
f'_2=\;\check{\kappa}_{p\dis\top}\cirk\hat{k}^{2}_{p,\bot}\;=\;
\check{k}^{2}_{p,\top}\cirk\hat{\kappa}_{p\kon\bot}\!.
\end{array}
\]
(This follows from Restricted Dicartesian Coherence too.) If we
write $\mO_{\bot,\top}$ for $\hat{\kappa}_{\bot}$, which is equal
to $\check{\kappa}_{\top}$ in $\ml_{\top,\bot}$, then in
$\ml_{\top,\bot}$ we have
\[
f'_2=\;\check{k}^{2}_{p,\top}\cirk
\mO_{\bot,\top}\cirk\hat{k}^{2}_{p,\bot}.
\]
So in $\ml_{\top,\bot}$ extended with $f_1=f_2$ we can derive
\[
(\hat{k}^{}\check{k}^{})\quad
\check{k}^{1}_{p,\top}\cirk\hat{k}^{1}_{p,\bot}\;=\;
\check{k}^{2}_{p,\top}\cirk
\mO_{\bot,\top}\cirk\hat{k}^{2}_{p,\bot}.
\]

The equation
\[
(\hat{k}^{}\check{\kappa})\quad
\hat{k}^{1}_{p,\bot}\;=\;\check{\kappa}_p\cirk\hat{k}^{2}_{p,\bot},
\]
which holds in $\set_{\ast}^{\emptyset}$, and which we have used
above for showing the nonmaximality of $\ml_{\top,\bot}$, clearly
yields $(\hat{k}^{}\check{k}^{})$, which hence holds in
$\set_{\ast}^{\emptyset}$, and which hence we could have also used
for showing this nonmaximality.

If we refine the procedure above by building $A'$ and $B'$ out of
$A$ and $B$ more carefully, then in some cases we could derive
$(\hat{k}^{}\check{\kappa})$ or its dual
\[
(\check{k}^{}\hat{\kappa})\quad
\check{k}^{1}_{p,\top}\;=\;\check{k}^{2}_{p,\top}\cirk\hat{\kappa}_p
\]
instead of $(\hat{k}^{}\check{k}^{})$. We do not replace $x$ by
$p\kon\bot$ in building $A'$, and we can proceed more selectively
with other occurrences of letters and $\top$ in $A$, in order to
obtain an $A'$ isomorphic to $p$ if possible.  We can proceed
analogously when we build $B'$ out of $B$ to obtain a $B'$
isomorphic to $p$ if possible.

Note that we have the following:
\[
\begin{array}{ll}
\check{\kappa}_{p\kon\bot}\cirk\hat{k}^{2}_{p,\bot}
\!\!\!\!&=\langle
\check{\kappa}_p,\mj_{\bot}\rangle\cirk\hat{k}^{2}_{p,\bot}
\\[.5ex]
&=\langle \hat{k}^{1}_{p,\bot},\hat{k}^{2}_{p,\bot}\rangle,
{\mbox{\rm{ with }}} (\hat{k}^{}\check{\kappa}),\\[.5ex]
&=\mj_{p\kon\bot}.
\end{array}
\]
In the other direction, it is clear that the equation derived
yields $(\hat{k}^{}\check{\kappa})$. So with
$(\hat{k}^{}\check{\kappa})$ we have that $C\kon\bot$ and $\bot$
are isomorphic, and, analogously, with
$(\check{k}^{}\hat{\kappa})$ we have that $C\dis\top$ and $\top$
are isomorphic. It can be shown that the natural logical category
defined as $\ml_{\top,\bot}$ save that we assume in addition both
$(\hat{k}^{}\check{\kappa})$ and $(\check{k}^{}\hat{\kappa})$ is
maximal. (This is achieved by eliminating letterless subformulae
from $C$ and $D$ in $g_1,g_2\!:C\vdash D$ such that $Gg_1\neq
Gg_2$, and falling upon the argument used for the maximality of
\ml\ in the preceding section.)

If $f\!:a\vdash b$ is any arrow of a dicartesian category \aA\ and
$(\hat{k}^{}\check{k}^{})$ holds in \aA, then we have in \aA
\[
\begin{array}{ll}
\check{k}^{1}_{b,\top}\cirk f\cirk\hat{k}^{1}_{a,\bot}
\!\!\!\!&=\;\check{k}^{1}_{b,\top}\cirk\hat{k}^{1}_{b,\bot}\cirk
(f\kon\mj_{\bot})
\\[.5ex]
&=\;\check{k}^{2}_{b,\top}\cirk
\mO_{\bot,\top}\cirk\hat{k}^{2}_{a,\bot},
\end{array}
\]
and hence for $f,g\!:a\vdash b$ we have in \aA
\[
(\hat{k}^{}\check{k}^{}fg)\quad \check{k}^{1}_{b,\top}\cirk
f\cirk\hat{k}^{1}_{a,\bot}\;=\; \check{k}^{1}_{b,\top}\cirk
g\cirk\hat{k}^{1}_{a,\bot}.
\]

So, although $\ml_{\top,\bot}$ is not maximal, it is maximal in
the relative sense that every dicartesian category that satisfies
an equation $f=g$ between arrow terms of $\ml_{\top,\bot}$ such
that $Gf\neq Gg$ satisfies also $(\hat{k}^{}\check{k}^{})$ and
$(\hat{k}^{}\check{k}^{}fg)$. Some of these dicartesian categories
may satisfy more than just $(\hat{k}^{}\check{k}^{})$ and
$(\hat{k}^{}\check{k}^{}fg)$. They may satisfy
$(\hat{k}^{}\check{\kappa})$ or $(\check{k}^{}\hat{\kappa})$,
which yields
\[
f\cirk\hat{k}^{1}_{a,\bot} \;=g\cirk\hat{k}^{1}_{a,\bot}
\quad\quad {\mbox{\rm{or }}} \quad\quad
\check{k}^{1}_{b,\top}\cirk f =\;\check{k}^{1}_{b,\top}\cirk g,
\]
and some may be preorders.

\vspace{2ex}

\noindent {\small {\it Acknowledgement$\,$}. We are grateful to
Slobodanka Jankovi\' c for a helpful stylistic suggestion. Work on
this paper was supported by the Ministry of Science of Serbia
(Grant 144013).}

\end{document}